\documentclass[journal]{IEEEtran}
\usepackage{amssymb}
\usepackage{amsmath}
\usepackage{amscd}
\usepackage{graphicx}
\usepackage{graphics}
\usepackage{xcolor,graphicx}
\usepackage{csquotes}
\usepackage{mathtools}
\usepackage{caption}
\usepackage{float}
\usepackage{subcaption}
 \usepackage{algorithmic}
 \usepackage{algorithm}
\usepackage{mathabx}
\ifCLASSINFOpdf
\else
\fi
\begin{document}
%
\title{Full Waveform Inversion with Adaptive Regularization}

%
%
%

\author{Hossein S. Aghamiry,~\IEEEmembership{}
        Ali Gholami,~\IEEEmembership{}
        St\'ephane Operto,~\IEEEmembership{}

\thanks{Hossein S. Aghamiry is with Universit\'e C\^ote d'Azur, CNRS, Observatoire de la C\^ote d'Azur, IRD , G\'eoazur, Valbonne, France,
  e-mail: aghamiry@geoazur.unice.fr. }
  \thanks{Ali Gholami is with the Institute of Geophysics, University of Tehran, Tehran, Iran,
  e-mail: agholami@ut.ac.ir.}
\thanks{St\'ephane Operto is with Universit\'e C\^ote d'Azur, CNRS, Observatoire de la C\^ote d'Azur, IRD , G\'eoazur, Valbonne, France,
  e-mail: operto@geoazur.unice.fr.}
         }
\maketitle

\begin{abstract}
Regularization is necessary for solving nonlinear ill-posed inverse problems arising in different fields of geosciences.
The base of a suitable regularization is the prior expressed by the regularizer, which can be non-adaptive or adaptive (data-driven). Nevertheless, tailoring a suitable and easy-to-implement prior for describing geophysical models is a nontrivial task.
In this paper, we propose general black-box regularization algorithms for solving nonlinear inverse problems such as full-waveform inversion (FWI), which admit \textit{empirical} priors that are determined \textit{adaptively} by sophisticated denoising algorithms. 
The nonlinear inverse problem is solved by a proximal Newton method, which generalizes the traditional Newton step in such a way to involve the gradients/subgradients of a (possibly non-differentiable) regularization function through operator splitting and proximal mappings.
Furthermore, it requires to account for the Hessian matrix in the regularized least-squares optimization problem. 
We propose two different splitting algorithms for this task. In the first, we compute the Newton search direction with an iterative method based upon the first-order generalized iterative shrinkage-thresholding algorithm (ISTA), and hence Newton-ISTA (NISTA). The iterations require only Hessian-vector products to compute the gradient step of the quadratic approximation of the nonlinear objective function. The second relies on the alternating direction method of multipliers (ADMM), and hence Newton-ADMM (NADMM), where the least-square optimization subproblem and the regularization subproblem in the composite are decoupled through auxiliary variable and solved in an alternating mode. The least-squares subproblem can be solved with exact, inexact, or quasi-Newton methods. We compare NISTA and NADMM numerically by solving full-waveform inversion with BM3D regularizations. The
tests show promising results obtained by both algorithms. However, NADMM shows a faster convergence rate than Newton-ISTA when using L-BFGS to solve the Newton system.
\end{abstract}
\graphicspath{{./Fig/}}
\section{Introduction}
\IEEEPARstart{N}ONLINEAR inverse problems frequently arise in different fields of geosciences \cite{Tarantola_2005_IPT}. Large-scale problems are typically solved with iterative local optimization (gradient-based) techniques such as Newton's method. Furthermore, such problems are inherently ill-posed and thus require regularization techniques to be implemented such that assumptions and priors about the unknown models are encoded in the optimization.
At the heart of a suitable regularization is a priori information expressed by the regularizer or regularization function \cite{Gholami_2010_RLN}. A proper regularizer, added to the objective function, renders the solution unique, increases its stability, and prevents data overfitting. It should be able to (mathematically) describe the solution while being easy to implement with iterative linearization methods.
These specifications make tailoring a suitable regularizer nontrivial.
A prior can be \textit{adaptive} or \textit{non-adaptive}, where by {\it{adaptive}} is meant the adaptation of the regularization function to the problem of interest.
Traditional priors used to solve inverse problems such as smoothness, sparseness, blockiness are \textit{non-adaptive} \cite{Tikhonov_1977_SIP,Tarantola_2005_IPT}. They are defined according to the preliminary assumptions about the targeted model, which are independent of the data and the problem to be solved. In contrast, \textit{adaptive} priors are solely derived from the data and tailored to the model accordingly.
Complex models require complex priors, which can be hard to derive.
Different priors lead to different forms of regularization, ranging from smooth and convex single-parameter regularizers \cite{Tikhonov_1977_SIP} to non-smooth and non-convex multi-parameter ones \cite{Gholami_2011_GFS,Selesnick_2017_SSA}. 

Denoising as the simplest inverse problem has contributed to enormous progress in developing sophisticated adaptive and non-adaptive priors for complicated signal recovery from noisy signals \cite{Milanfar_2012_TMI}. 
Some recently proposed excellent denoising methods include nonlocal means filters \cite{Milanfar_2012_TMI,Goyal_2020_IDR} and block matching 3D filter (BM3D) \cite{Dabov_2007_VDB} and its variants \cite{Goyal_2020_IDR}.
These patch-based methods use both local and nonlocal redundancy of information in the input signal to preserve structures in the solution by yielding locally adaptive filters via similarity kernels. Specifying the kernel function in these methods is essentially equivalent to estimating a particular type of empirical prior from the input signal \cite{Milanfar_2012_TMI}.
This somehow contrasts with the traditional non-adaptive regularization methods, for which the prior is fixed and independent from the observed data \cite{Tarantola_2005_IPT}.
Such an adaptive regularization has been applied to linear inverse problems in, e.g., \cite{Danielyan_2011_BFV} and \cite{Venkatakrishnan_2013_PAP}.

In this paper, we extend such adaptive methods to nonlinear inverse problems via adaptive proximal Newton-type algorithms.  
Similar to the classical Newton-type methods, a nonlinear inverse problem is solved iteratively as 
$\bold{m}_{k+1}= \bold{m}_{k}+\alpha_k \Delta\bold{m}_{k}$, where $\bold{m}_{k}$ is the model parameters at iteration $k$, $ \Delta\bold{m}_{k}$ is the search direction and $\alpha_k$ is the step length.
%
When a composite objective function includes a general (and possibly non-differentiable) regularization term,
$\Delta\bold{m}_{k}$ must further involve the gradients/subgradients of the regularization function \cite{Lee_2014_PDG}. Proximal Newton methods achieve this task by breaking down the original complex problem into simpler subproblems through operator splitting and proximal mappings.
We propose two distinct algorithms to solve the regularized problem with proximal Newton methods.
In the first, called NISTA, the Newton search direction $\Delta\bold{m}_{k}$ at iteration $k$ is computed by minimizing a composite objective function given by the sum of the locally quadratic approximation of the nonlinear misfit function involving the Hessian and the regularization function. The minimum of this surrogate objective function is found iteratively with a proximal gradient method based upon the shrinkage-thresholding algorithm (ISTA) \cite{Daubechies_2004_ITA, Attouch_2016_RCN}. 
A key property of this method is to require only the Hessian-vector product to build the gradient of the linearized misfit function. 
  
The second algorithm, called NADMM, relies on the alternating direction method of multipliers (ADMM) \cite{Boyd_2011_DOS,Aghamiry_2019_IWR}.
ADMM decouples the linearized least-squares objective function and the regularization term via an auxiliary variable and solves the two subproblems in alternating mode with the primal-dual method of multipliers.
The first subproblem requires to solve a linear system involving the Hessian, just like classical Newton-type methods. This system can be solved exactly or approximately with inexact or quasi-Newton algorithms.
%
%

An important property of the proposed algorithms is that they only need the outputs of the regularizer without asking for any information about its functional form and statistical properties. This black-box implementation brings flexibility to the methods for using adaptive and non-adaptive regularizations with the local solvers.
The main properties of the proposed regularization are as follows: [1] It can be easily implemented with existing algorithms. [2] The regularizer is treated as black-box, and thus, adaptive and non-adaptive regularizations can be implemented. [3] Irrespective of the differentiability of the regularizer, it can be implemented with the iterative gradient-based solvers. [4] The computational overhead generated by the regularization is the computation of the proximal/denoising operator at each iteration, and hence is negligible in most cases. 

We implement the proposed adaptive regularization to solve full-waveform inversion (FWI), an ill-posed PDE-constrained nonlinear optimization problem, in which the subsurface parameters and the wavefields are defined as the minimizers of the Euclidean distance between observed and calculated data \cite{Tarantola_1984_ISR,Pratt_1998_GNF}. Among different methods to solve this constrained optimization problem, we consider a variable projection formulation leading to the classical FWI \cite{Pratt_1998_GNF} and ADMM \cite{Aghamiry_2019_IWR}. 
The ADMM formulation, which updates the parameters and the wavefields in alternating mode, is referred to iteratively refined wavefield reconstruction inversion (IR-WRI) \cite{Aghamiry_2019_IWR}.  
Numerical tests performed show outstanding performance of the adaptive regularization in building complicated velocity models by the above waveform inversion methods. 

%
\section{Preliminaries}
As in this paper we will use the concepts and formulas used in linear inverse theory. A brief review of these concepts is given here.

In linear inverse problems the desired model, denoted by column vector $\bold{m}$, needs to be estimated from measurements $\bold{d}$ that relate to $\bold{m}$ via a linear operator/matrix $\bold{A}$, i.e. $\bold{d}=\bold{Am+e}$ for some random noise $\bold{e}$.  
For a Gaussian distributed random noise, the estimation problem usually appears as determination of the minimizer of a suitably defined objective function 
\begin{equation} \label{lin}
\arg\min_{\bold{m}} \frac12 \|\bold{d}-\bold{Am}\|_2^2 + \lambda \mathcal{R}(\bold{m}),
\end{equation}
where $\mathcal{R}$  is a regularizer or regularization function which somehow prevents data overfitting and $\lambda$ determines  regularization weight. Different forms of $\mathcal{R}$ have been proposed, ranging from smooth and convex single-parameter functions \cite{Tikhonov_1977_SIP} to non-smooth and non-convex multi-parameter ones \cite{Gholami_2011_GFS,Selesnick_2017_SSA}.
In its simplest form $\mathcal{R}(\bold{m})=\|\bold{m}-\bold{m}^{prior}\|_2^2$ is a damping term that encourages $\bold{m}$ not to be very far from the prior model $\bold{m}^{prior}$ \cite{Tarantola_2005_IPT}.  
\subsection{Denoising and Proximal Operator}
In denoising problem $\bold{A=I}$ (the identity matrix) and the estimate is simply  defined as
\begin{equation} \label{prox0}
\text{prox}_{\lambda \mathcal{R}}(\bold{d})=\arg\min_{\bold{m}} \frac12 \|\bold{d}-\bold{m}\|_2^2 + \lambda \mathcal{R}(\bold{m}).
\end{equation}
This is called the proximal operator of $\mathcal{R}$ \cite{Combettes_2011_PRO}.
Despite its simple definition, proximal operators are powerful tools in optimization because 
1) the general optimization problem \eqref{lin} can be solved by proximal algorithms which merely require to evaluate the gradient of the misfit function, $\mathcal{M}(\bold{m})=\frac12\|\bold{d}-\bold{Am}\|_2^2$, and the proximal operator \eqref{prox0}. 2) Since a proximal operator involves the information about gradients/subgradients of $\mathcal{R}$, proximal algorithms handle both differentiable and nondifferentiable forms of $\mathcal{R}$. This is in contrast with the traditional algorithms, such as Newton's algorithm, which requires the objective to be differentiable.
Furthermore, the interpretation of the proximal operator as a denoising \cite{Kamilov_2017_PPA} allows us to solve  \eqref{lin} with advanced regularizations embedded in sophisticated denoising algorithms such as non-local means (NLM), block-matching 3D (BM3D) or deep learning denoisers \cite{Meinhardt_2017_LPO}.

\subsection{The Proximal Gradient Method}
The proximal-gradient method is an important tool for solving non-linear problems we describe in subsequent sections.
In order to see how the proximal operator \eqref{prox0} helps to solve \eqref{lin}, we use the majorization-minimization (MM) approach \cite{Lange_2016_MOA} which has a simple convergence proof. It is interesting to note that for $\mathcal{M}(\bold{m})=\frac12\|\bold{d}-\bold{Am}\|_2^2$ and $c \in (0,1/\vvvert\bold{A}\vvvert^2)$, with $\vvvert\bold{A}\vvvert$ the largest singular value of $\bold{A}$, we have that
\begin{equation} \label{lin_majo}
\mathcal{M}(\bold{m}) +\lambda \mathcal{R}(\bold{m}) \leq \widetilde{\mathcal{M}}_k(\bold{m})+\lambda \mathcal{R}(\bold{m})
\end{equation}
with equality at $\bold{m}=\bold{m}_k$, where
\begin{equation} \label{majo}
\widetilde{\mathcal{M}}_k(\bold{m})= \mathcal{M}(\bold{m}_k)+(\bold{m}-\bold{m}_k)^T\nabla\mathcal{M}(\bold{m}_k)
+\frac{1}{2c}\|\bold{m}-\bold{m}_k\|_2^2, 
\end{equation}
in which $\bold{m}_k$ is a reference model (previous iterate) and $\nabla\mathcal{M}(\bold{m}_k)$ is the gradient vector.
This approximation allows us to minimize \eqref{lin} by iteratively minimizing a simpler problem
\begin{equation} \label{lin_sub}
\bold{m}_{k+1}=\arg\min_{\bold{m}} \widetilde{\mathcal{M}}_k(\bold{m}) + \lambda \mathcal{R}(\bold{m}).
\end{equation}
Simple algebra shows that \eqref{lin_sub} is equivalent to
\begin{equation} \label{ISTA}
\bold{m}_{k+1}=\text{prox}_{c\lambda \mathcal{R}}(\bold{m}_{k}- c \nabla\widetilde{\mathcal{M}}_k(\bold{m}_k)).
\end{equation}
This is nothing but the famous iterative shrinkage-thresholding algorithm (ISTA) \cite{Daubechies_2004_ITA} (also known as
forward-backward splitting algorithm and proximal gradient method).
FISTA \cite{Beck_2009_FIS} is an accelerated version that uses a particular linear combination of the two last iterates to perform the update.
A simple acceleration is obtained by using the extrapolation method of Nesterov \cite{Nesterov_1983_AMO}, leading to the generalized form of ISTA \cite{Attouch_2016_RCN}
\begin{equation} \label{GISTA}
\begin{cases}
\bold{m}_{k+1}=\text{prox}_{c\lambda \mathcal{R}}(\bold{p}_{k}- c \nabla\widetilde{\mathcal{M}}_k(\bold{p}_k))  \\
\bold{p}_{k+1} =\bold{m}_{k+1} + \frac{k-1}{k+2} (\bold{m}_{k+1} - \bold{m}_{k}).
\end{cases}
\end{equation}
\section{Method} 
A nonlinear inverse problem such as FWI with a general form of regularization can be written as
\begin{equation} \label{main_gen}
\min_{\bold{m}} ~ \mathcal{M}(\bold{m}) + \lambda \mathcal{R}(\bold{m}),
\end{equation}
where $\bold{m}$ is the model parameters.
In \eqref{main_gen}, $\mathcal{M}(\bold{m})$ is the data misfit function. Its minimization ensures that the simulated data $F(\bold{m})$ are close to the measurements $\bold{d}$, where $F$ is a nonlinear differentiable function. 
$\mathcal{R}(\bold{m})$ is the possibly non-differentiable regularization, which encodes the prior knowledge about the model parameters and prevents data overfitting.
 $\lambda$ is the trade-off parameter that balances between the data misfit and regularization terms.

A Newton-type method approximates the misfit term with a local quadratic function of form
\begin{eqnarray} \label{quad_app}
&&{\widetilde{\mathcal{M}}_k}(\bold{m}) = \mathcal{M}(\bold{m}_k) + (\bold{m}-\bold{m}_k)^T\nabla\mathcal{M}(\bold{m}_k) \nonumber \\
&&\hspace{2.5cm}+\frac12 (\bold{m}-\bold{m}_k)^T\bold{H}_k(\bold{m}-\bold{m}_k), 
\end{eqnarray}
where $\bold{m}_k$ is the iterate at iteration $k$, $\nabla\mathcal{M}(\bold{m}_k)$ is the gradient vector, and $\bold{H}_k$ is the Hessian matrix $\nabla^2\mathcal{M}(\bold{m}_k)$ or an approximation of it. 

Using the approximation in \eqref{quad_app}, proximal Newton-type methods solve problem \eqref{main_gen} iteratively as
\begin{equation} \label{PNOPT}
\bold{m}_{k+1} = \bold{m}_{k} + \alpha_k \Delta\bold{m}_{k},
\end{equation}
where $\alpha_k$ is the step length, which an be determined by a line search method, and
\begin{align} \label{PNOPT_g}
\Delta\bold{m}_{k}  &= \arg\min_{\Delta\bold{m}} ~\widetilde{\mathcal{M}}_k(\bold{m}_k+\Delta\bold{m}) + \lambda \mathcal{R}(\bold{m}_k+\Delta\bold{m})
\end{align} 
is a search direction \cite{Lee_2014_PDG}.
Computation of the search direction $\Delta\bold{m}_{k}$ is the most computationally expensive part of this algorithm because it requires the minimization of a composite function given by the sum of a quadratic term involving the Hessian matrix, \eqref{quad_app}, and the regularization term $\mathcal{R}$. For $\lambda=0$, the algorithm reduces to a classical Newton method, where an approximation of the Hessian can be employed, leading to quasi-Newton methods or gradient method if $\bold{H}_k$ reduces to a scaled version of the identity matrix, \eqref{majo}.
%
%
For $\lambda\neq 0$, however, determination of the search direction in \eqref{PNOPT_g} is more challenging. 
In the following, we propose two methods for this task.

\subsection{Newton-ISTA (NISTA)}
NISTA relies on the first-order ISTA, \eqref{ISTA}, to estimate iteratively the Newton search direction \eqref{PNOPT_g}. This requires to implement the following inner loop  within the outer loop over $k$ 
\begin{align}
\label{nista}
\begin{cases}
\Delta\bold{m}_k^{\ell+\frac12} = \Delta\bold{p}^{\ell} - c_k (\bold{H}_k\Delta\bold{p}^{\ell} + \nabla\mathcal{M}(\bold{m}_k))\\
\Delta\bold{m}_k^{\ell+1} =\text{prox}_{c_k\lambda \mathcal{R}}(\bold{m}_{k}+\Delta\bold{m}_k^{\ell+\frac12}) - \bold{m}_{k}\\
\Delta\bold{p}^{\ell+1} =\Delta\bold{m}_{k}^{\ell+1} + \frac{\ell-1}{\ell+2} (\Delta\bold{m}_{k}^{\ell+1} - \Delta\bold{m}_{k}^{\ell}),
\end{cases}
\end{align}
where $l$ is the inner iteration count, $\Delta\bold{p}^{0}=0$, and $c_k \in (0,1/\vvvert\bold{H}_k\vvvert^2)$. The term in bracket in the first line of \eqref{nista} is the gradient of the surrogate function $\widetilde{\mathcal{M}}_k(\bold{m})$, \eqref{quad_app}.
For many choices of the regularizer $\mathcal{R}$, there can be a closed-form expression for the denoiser in the second subproblem of \eqref{nista}.
The main property of this formulation is that it can be generalized to exploit multiple (even data-driven)  priors by using different denoisers instead of the prox operator, e.g. BM3D \cite{Dabov_2007_VDB}.
It is important that the denoiser function is treated as a black box, i.e., we only need access to the output of the denoiser for a given input, irrespective of its functional form. 

The NISTA is summarized in Algorithm \ref{alg1}. The algorithm is started with $\Delta\bold{p}_0=\bold{0}$. However, to improve the convergent speed, we can perform a warm start of the inner loop by using the results of the previous iteration.

\begin{algorithm}
 \begin{algorithmic}[1]
 \caption{Adaptive regularization by NISTA.} \label{alg1}
 \REQUIRE starting point $\bold{m}_k$
 \STATE set $\Delta\bold{p}^0=\bold{0}$
 \REPEAT
 \STATE Compute the Hessian $\bold{H}_k$ or an approximation to it.
 \STATE Compute the step direction:
\FOR{$\ell = 1$ to $N-1$}
\STATE $\Delta\bold{m}_k^{\ell+\frac12} = \Delta\bold{p}^{\ell} - c_k (\bold{H}_k\Delta\bold{p}^{\ell} + \nabla\mathcal{M}(\bold{m}_k))$
\STATE $\Delta\bold{m}_k^{\ell+1} =\text{prox}_{c_k\lambda \mathcal{R}}(\bold{m}_{k}+\Delta\bold{m}_k^{\ell+\frac12}) - \bold{m}_{k}$
\STATE $\Delta\bold{p}^{\ell+1} =\Delta\bold{m}_{k}^{\ell+1} + \frac{\ell-1}{\ell+2} (\Delta\bold{m}_{k}^{\ell+1} - \Delta\bold{m}_{k}^{\ell})$
\ENDFOR
\STATE Select step length $\alpha_k$ with a backtracking line search.
\STATE Update: $\bold{m}_{k+1} = \bold{m}_{k} + \alpha_k \Delta\bold{m}_k^N$.
\UNTIL{stopping conditions are satisfied.}
\end{algorithmic}
\end{algorithm}

\subsection{Newton-ADMM (NADMM)}
NADMM is obtained by solving \eqref{PNOPT_g} via the alternating direction method of multipliers (ADMM) \cite{Boyd_2011_DOS}. By introducing the auxiliary variable $\bold{p}=\bold{m}_k+\Delta\bold{m}$, we recast the minimization problem in \eqref{PNOPT_g} as the following constrained problem:
\begin{align} \label{main_eq}
&\min_{\Delta\bold{m,p}} ~~~~~~~\widetilde{\mathcal{M}}_k(\bold{m}_k+\Delta\bold{m}) + \lambda \mathcal{R}(\bold{p})\\
&\text{subject to}~~~ \bold{m}_k+\Delta\bold{m}=\bold{p}. \nonumber
\end{align}
Solving \eqref{main_eq} with an augmented Lagrangian method leads to the following saddle point problem
\begin{align} \label{main_al}
\min_{\Delta\bold{m,p}} \max_{\bold{q}}~&\widetilde{\mathcal{M}}_k(\bold{m}_k+\Delta\bold{m}) + \lambda \mathcal{R}(\bold{p})\\
&+ \langle \bold{q},\bold{m}_k+\Delta\bold{m}-\bold{p}\rangle
+ \frac{1}{2c_k}\|\bold{m}_k+\Delta\bold{m}-\bold{p}\|_2^2, \nonumber
\end{align}
where $\bold{q}$ is the Lagrange multiplier and $1/c_k$ serves as a penalty parameter. Applying the scaled form of ADMM to \eqref{main_al}, when combined with \eqref{PNOPT}, gives the iteration
\begin{align} 
\label{ADMM}
\begin{cases}
\Delta\bold{m}_k=\arg\min_{\Delta\bold{m}}~\widetilde{\mathcal{M}}_k(\bold{m}_k+\Delta\bold{m}) \\
 \hspace{2.75cm}+ \frac{1}{2c_k}\|\bold{m}_k+\Delta\bold{m}-\bold{p}_k-\bold{q}_k\|_2^2 \\
\bold{m}_{k+1} = \bold{m}_{k} + \alpha_k \Delta \bold{m}_{k} \\
\bold{p}_{k+1}=\text{prox}_{c_k\lambda \mathcal{R}}(\bold{m}_{k+1}-\bold{q}_k)\\
\bold{q}_{k+1}=\bold{q}_k+\bold{p}_{k+1} - \bold{m}_{k+1},
\end{cases}
\end{align}
where the primal and dual variables are updated in alternating mode.
With a change of variable $\bold{m}_k^{prior}=\bold{p}_k+\bold{q}_k$, the first subproblem in \eqref{ADMM} requires us to solve 
\begin{equation} \label{M_sub}
\min_{\Delta\bold{m}} ~\widetilde{\mathcal{M}}_k(\bold{m}_k+\Delta\bold{m}) + \frac{1}{2c_k} \|\bold{m}_k+\Delta\bold{m}-\bold{m}_k^{prior}\|_2^2
\end{equation}
which has a closed-form minimizer given by
\begin{equation}
\Delta\bold{m}_{k} = (c_k\bold{H}_k+ \bold{I})^{-1}(-c_k\nabla\mathcal{M}(\bold{m}_k) + \Delta\bold{m}_k^{prior}),
\label{eqsubadmm}
\end{equation} 
where $\Delta\bold{m}_k^{prior}=\bold{m}_k^{prior}-\bold{m}_k$.
This is a {\it generalized} gradient step because it implicitly includes the information carried out by the gradient/subgradient of the possibly non-differentiable regularizer. 
It is seen that the priori information introduced by $\mathcal{R}(\bold{m})$ in the original problem \eqref{main_gen}, regardless of its mathematical form or its differentiability, is replaced by a priori information that the (unknown) model at each iteration is a sample of a known Gaussian probability density whose mean is $\bold{m}_k^{prior}$ and whose covariance matrix is a scaled identity matrix. The regularization appeared as a damping term that encourages the model not to be very far from the dynamic prior/reference model $\bold{m}_k^{prior}$, unlike traditional Bayesian approach \cite{Tarantola_2005_IPT} where the a priori model is static. The Newton system \eqref{eqsubadmm} can be solved with any quasi-Newton or inexact Newton methods.

The proposed NADMM method is summarized in Algorithm \ref{alg2}.

\begin{algorithm}
 \begin{algorithmic}[1]
 \caption{Adaptive regularization by NADMM.} \label{alg2}
 \REQUIRE starting point $\bold{m}_0$
 \STATE set $\bold{p}_0=\bold{q}_0=\bold{0}$
 \REPEAT
 \STATE Compute the Hessian $\bold{H}_k$ or an approximation to it.
 \STATE Compute the step direction:\\
        $\Delta \bold{m}_{k} = (c_k\bold{H}_k+\bold{I})^{-1}(-c_k\nabla\mathcal{M}(\bold{m}_k) + \bold{p}_k+\bold{q}_k-\bold{m}_k).$
\STATE Select step length $\alpha_k$ with a backtracking line search.
\STATE Update: $\bold{m}_{k+1} = \bold{m}_{k} + \alpha_k \Delta \bold{m}_{k}$.
\STATE Update: $\bold{p}_{k+1} = \text{prox}_{c_k\lambda \mathcal{R}}(\bold{m}_{k+1}-\bold{q}_k).$
\STATE Update: $\bold{q}_{k+1}= \bold{q}_k + \bold{p}_{k+1}-\bold{m}_{k+1}$
\UNTIL{stopping conditions are satisfied.}
\end{algorithmic}
\end{algorithm}

\subsection{Application to Full Waveform Inversion}

In the {\it{Numerical example}} section, we assess the algorithms 1 and 2 against seismic full waveform inversion methods with a series of benchmarks. Here, we briefly review the two different formulations of full waveform inversion that will be used. The first classical one relies on variable projection to recast the nonlinear constrained problem as an unconstrained problem with a reduced search space. The second extends the linear regime of the waveform inversion with ADMM. This recasts the original nonlinear constrained problem as a biconvex problem according to the bilinearity of the wave equation.

\subsubsection{Reduced-space FWI}
In classical full waveform inversion \cite{Pratt_1998_GNF}
\begin{equation} \label{mis_fwi}
\mathcal{M}(\bold{m}) = \frac12 \|\bold{d}-F(\bold{m})\|_2^2,
\end{equation}
where $\bold{d}$ is the observed data and $F(\bold{m})=\bold{PA}^{-1}(\bold{m})\bold{b}$ is the calculated data in which $\bold{P}$ is the observation operator that samples the wavefield $\bold{A}^{-1}(\bold{m})\bold{b}$ at the receiver positions, $\bold{b}$ is the source and $\bold{A}$ is the wave-equation operator.

For \eqref{mis_fwi}, the gradient and the Hessian are given by \cite{Pratt_1998_GNF} 
\begin{eqnarray}
\nabla\mathcal{M}(\bold{m})= -\bold{J}^T\bold{\Delta d}
\end{eqnarray}
and
\begin{eqnarray} \label{FWI_Hess}
\nabla^2\mathcal{M}(\bold{m})= \bold{J}^T\bold{J} + \frac{\partial \bold{J}^T}{\partial \bold{m}^T},
[\bold{\Delta d}|\cdots |\bold{\Delta d}],
\end{eqnarray}
where $\bold{\Delta d}=\bold{d}-F(\bold{m})$ and $\bold{J}$ is the sensitivity or the Fr\'{e}chet derivative matrix, defined as
\begin{equation}
\bold{J}_{ij} = \frac{\partial [F(\bold{m})]_i}{\partial \bold{m}_j}.
\end{equation}


\subsubsection{ADMM-based Wavefield Reconstruction Inversion (IR-WRI)}
In classical FWI, the wave-equation $\bold{A}(\bold{m})\bold{u}=\bold{b}$ is solved exactly at each iteration to generate the reduced form of the objective function \eqref{mis_fwi}. 
In the wavefield reconstruction inversion (WRI) method, the wave-equation is satisfied approximately through a penalty method such that the simulated wavefields match the observations. Then, the parameters are updated from the wavefields by least-squares minimization of the wave equation errors \cite{VanLeeuwen_2013_MLM,Aghamiry_2019_IWR}. Updating the wavefields and the subsurface parameters in alternating mode at iteration $k$  leads to the following objective function for $\bold{m}$
\begin{equation} \label{mis_wri}
\mathcal{M}(\bold{m}) = \frac{1}{2} \|\bold{b}-\bold{A}(\bold{m})\bold{u}_k\|_2^2,
\end{equation}
where the so-called data-assimilated wavefield $\bold{u}_k$ is the least-squares solution of the overdetermined system gathering the wave equation and the observation equation
\begin{equation}
\begin{pmatrix}
\bold{A}(\bold{m}_{k-1})\\
\mu \bold{P}
\end{pmatrix}   \bold{u}_k =
\begin{pmatrix}
\bold{b}\\
\mu \bold{d}
\end{pmatrix},
\label{eqas}
\end{equation}
where 
$\mu>0$ is the penalty parameter. Note that \eqref{mis_wri} and \eqref{eqas} are provided assuming a single source experiment. For multiple sources, the objective function \eqref{mis_wri} is simply obtained by summation over sources,  while one augmented system \eqref{eqas} per source needs to be solved. 
For \eqref{mis_wri}, the gradient and Hessian are given by
\begin{eqnarray}
\nabla\mathcal{M}(\bold{m})= - \bold{L}^T(\bold{b}-\bold{A}(\bold{m})\bold{u}_{k})
\end{eqnarray}
and
\begin{eqnarray} \label{WRI_Hess}
\nabla^2\mathcal{M}(\bold{m})= \bold{L}^T\bold{L},
\end{eqnarray}
where
\begin{equation}
\bold{L}_{ij} = \frac{\partial [\bold{A}(\bold{m})\bold{u}_{k}]_i}{\partial \bold{m}_j}.
\end{equation}
\section{Numerical examples}
\subsection{A Toy Example}
We first show the performance of the proposed algorithm with a simple two-dimensional nonlinear optimization problem.
\begin{equation} \label{obj2d}
\min_{m_1,m_2} ~75(m_2-m_1^2)^2+(1-m_1)^2 + \lambda(|m_1|+|m_2|).
\end{equation}
Comparing this objective function with the canonical form \eqref{main_gen}, we get that $\mathcal{M}$ is the Rosenbrock function
\begin{equation} \label{ros}
\mathcal{M}(m_1,m_2)=75(m_2-m_1^2)^2+(1-m_1)^2, 
\end{equation}
and $\mathcal{R}$ is the $l_1$-norm
\begin{equation}
\mathcal{R}(m_1,m_2)=|m_1|+|m_2|.
\end{equation}
The Rosenbrock function is continuously differentiable and has a global minimum at (1,1). Adding the sparsity-promoting regularization term to this function however moves this global minimum toward zero in an specific path. 
For $0\leq \lambda \leq 3/2$ the global minimum occurs at $(m_1^*,m_2^*)$ where $m_1^*=(2-\lambda)/(2+2\lambda)$ and $m_2^*=(m_1^*)^2-\lambda/150$. For $3/2< \lambda \leq 2$ it occurs at $(m_1^*,0)$ where $m_1^*$ solves $300(m_1^*)^3+2m_1^*+\lambda-2=0$ and for $\lambda >2$ the function reaches its global minimum at (0,0).

\begin{figure}[!htb]
\center
\includegraphics[width=0.49\textwidth]{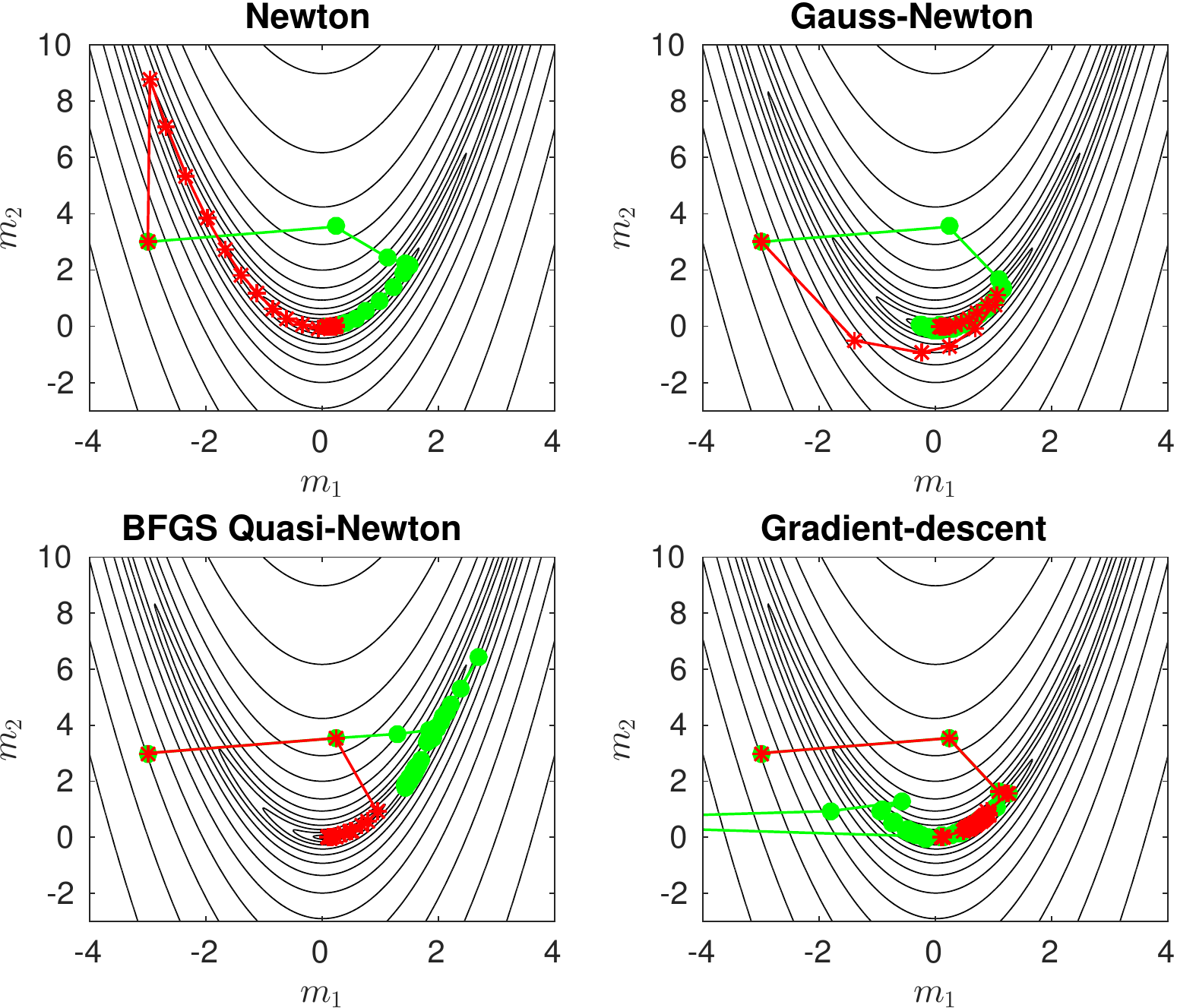}
\caption{Minimization of the sparsity-promoting regularized Rosenbrock function in \eqref{obj2d} via NISTA (green squares) and NADMM (red circles) using different approximations to the Hessian.}
\label{fig:Rosenbrock1}
\end{figure}

\begin{figure} [!htb]
\center
\includegraphics[width=0.49\textwidth]{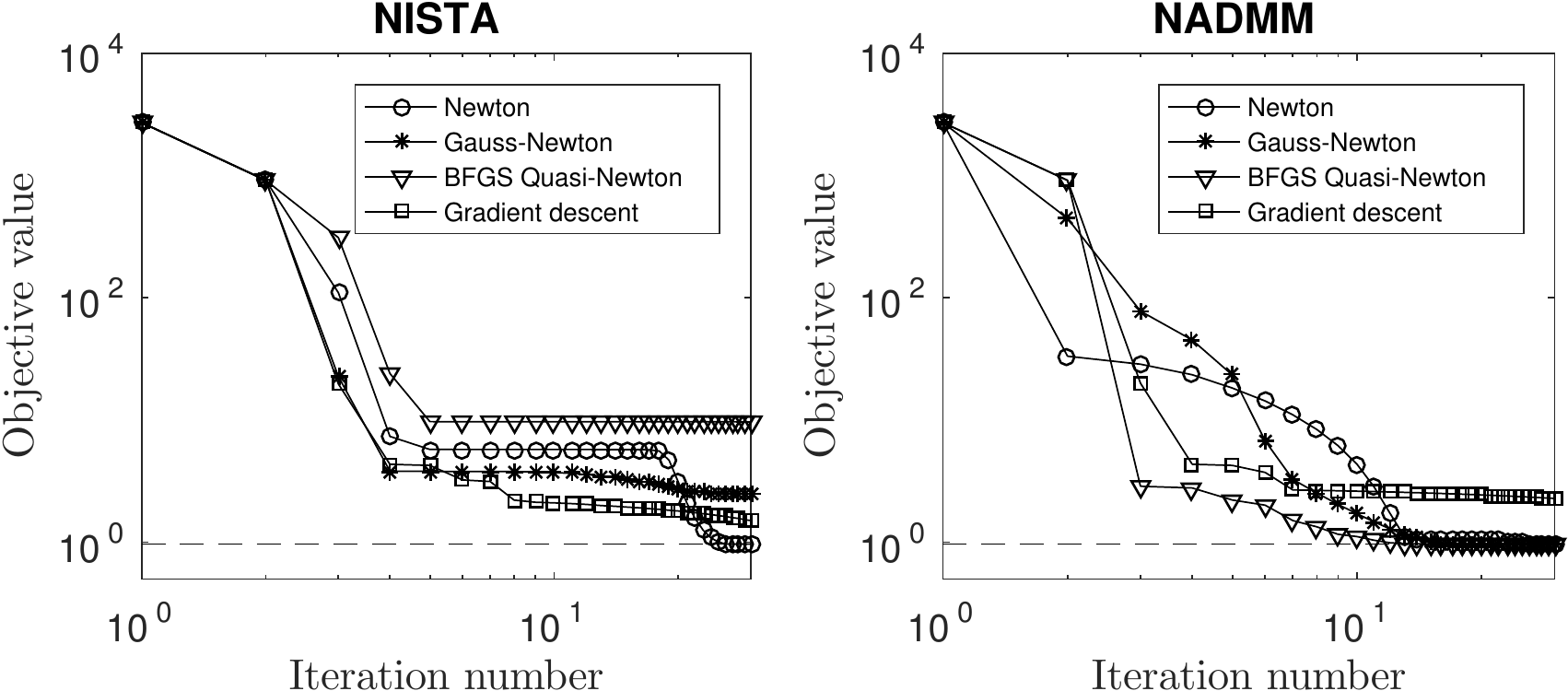}
\caption{Comparison between different methods in minimization of the sparsity regularized Rosenbrock function in \eqref{obj2d} via NISTA and NADMM.}
\label{fig:Rosenbrock2}
\end{figure}

We applied both NISTA and NADMM to minimize this nonlinear and nondifferentiable function for $\lambda=3/2$ (having the global minimum at (0.1,0), Fig. \ref{fig:Rosenbrock1}). We also use different approximations of the Hessian in each algorithm. 
The performance of all methods is compared in Fig. \ref{fig:Rosenbrock2} and as seen from this figure all methods successfully converged to the desired global minimum but with different number of iterations.
For all Newton-type methods, NADMM converges faster than NISTA (with 50 inner iterations).

\begin{figure*}[!htb]
\center
\includegraphics[width=12cm,clip=true,trim=0cm 0cm 0cm 0cm]{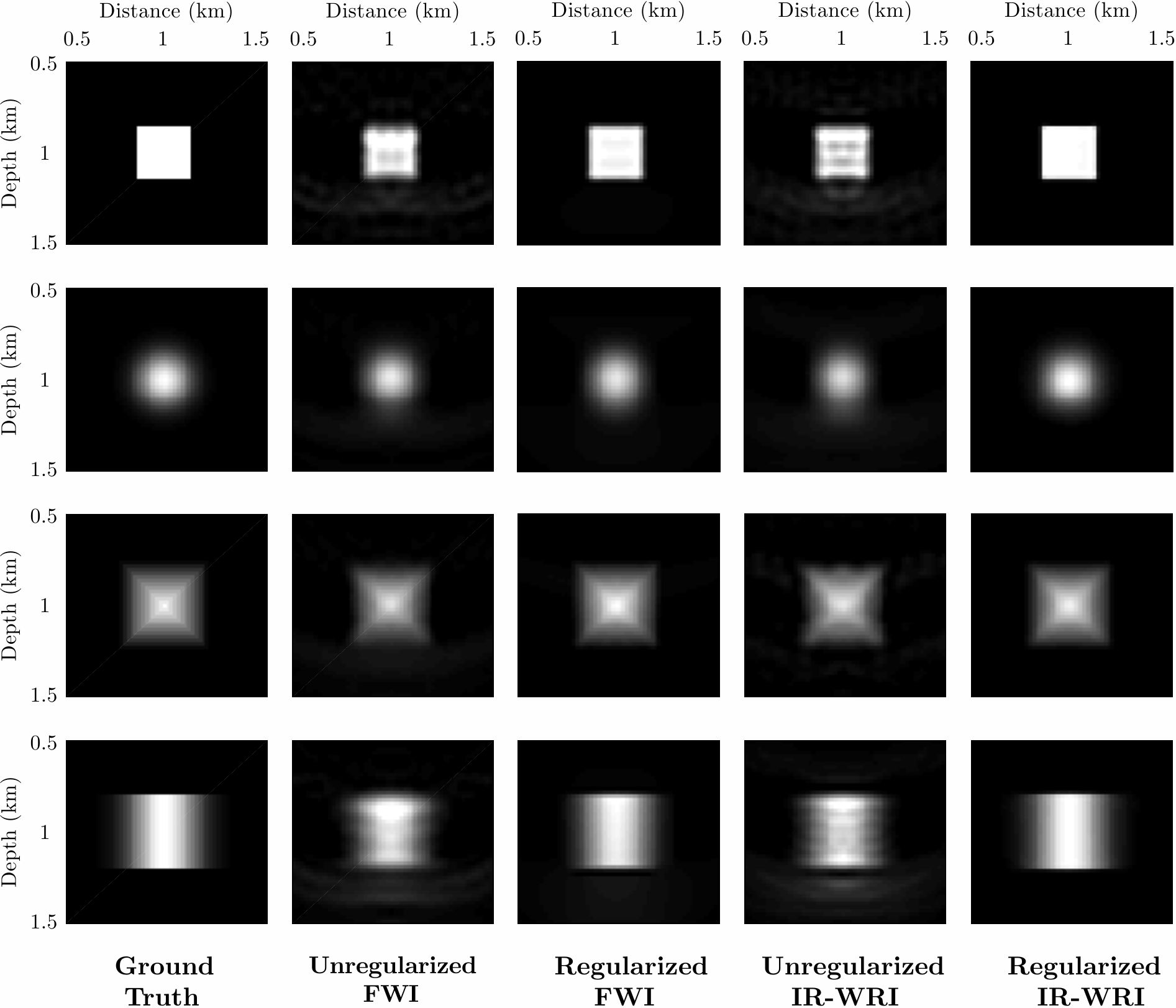}
\caption{The performances of FWI and IR-WRI with adaptive NADMM for reconstruction of different velocity structures. The data were generated by five sources at the surface (with 400~m spacing) and 50~m equally spaced receivers positioned on all the boundaries except the surface. In all figures, the colorbar varies between 2000 m/s and 2500 m/s with  low velocities in black and high velocities in white.}
\label{fig:Inclusions}
\end{figure*}

\subsection{Wavefield Inversion of Inclusions Models}
\begin{figure}[!htb] 
\includegraphics[width=0.49\textwidth]{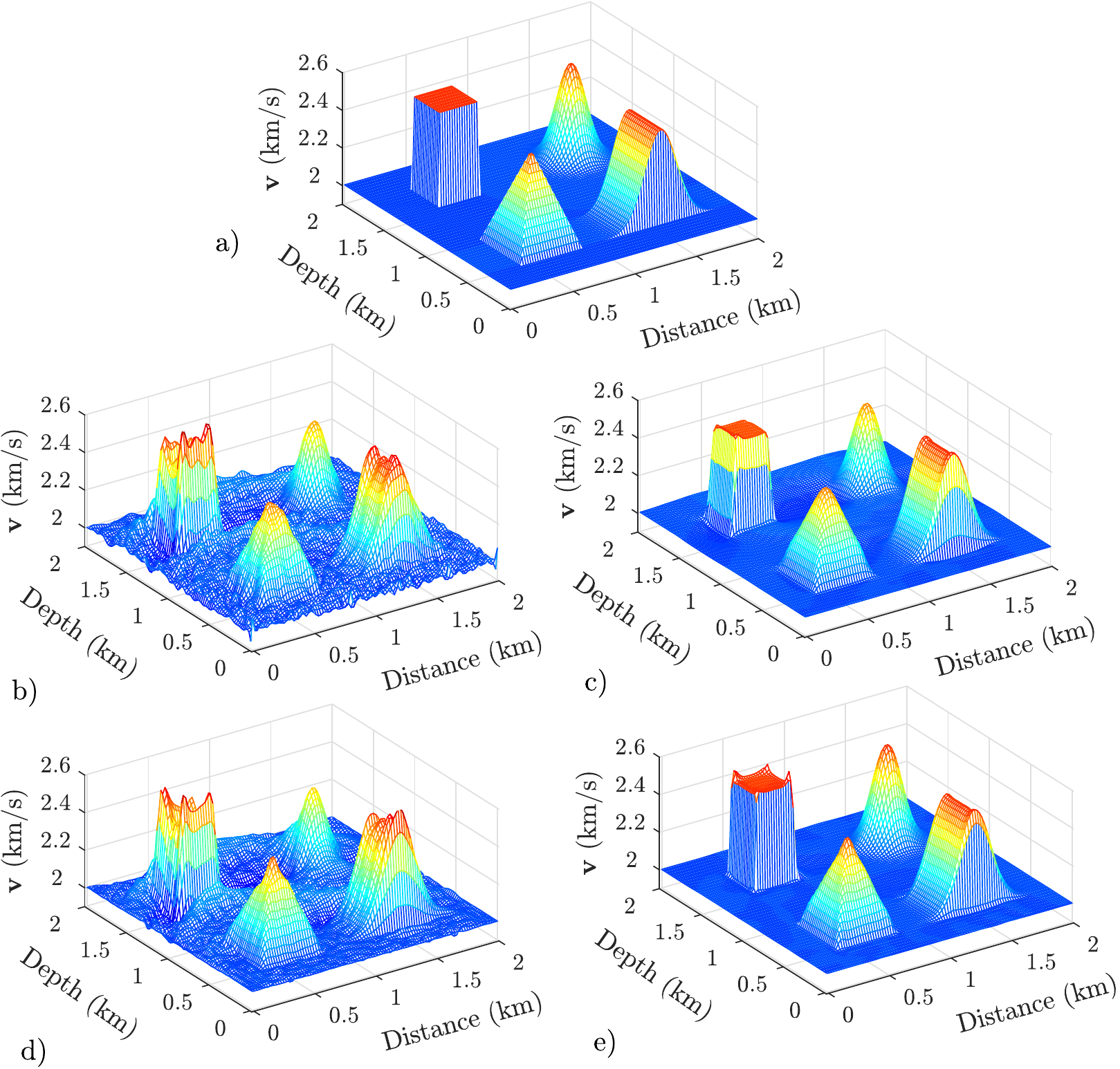}\\
\caption{Inclusion test. (a) True velocity model. (b-e) Velocity models estimated by (b) FWI without regularization, (c) FWI with regularization, (d) IR-WRI without regularization, (e) IR-WRI with regularization.}
\label{fig:simple}
\end{figure}

We now show how the proposed adaptive regularization helps us to construct different velocity models via FWI and IR-WRI when BM3D is used as denoiser. For this, we use four different velocity models, where the subsurface is  2 km $\times$ 2 km homogeneous velocity model ($\bold{v}=2$~km/s) including different inclusions with different characteristics (first column of Fig. \ref{fig:Inclusions}). Also, to show the flexibility of the proposed adaptive regularization in managing different priors simultaneously, we put all the four inclusions together in a model (Fig. \ref{fig:simple}a). For all the tests of this section, data are generated by five sources at the surface (with 400~m spacing) and 50~m equally spaced receivers placed on all the boundaries except the surface. 
%
%
The forward modeling is performed with a 9-point stencil finite-difference method implemented with anti-lumped mass and PML absorbing boundary conditions to solve the Helmholtz equation, where the stencil coefficients are optimized to the frequency \cite{Chen_2013_OFD} (this scenario is considered for all wave-propagation examples in this paper). The source signature is a Ricker wavelet with a 10~Hz dominant frequency. We start the inversion from the homogeneous background model ($\bold{v}=2$~km/s) and invert simultaneously four frequency components (5, 7, 10, and 12.5 Hz) with noiseless and noisy data for FWI and IR-WRI when the new regularization is used or not. 

We first apply FWI and IR-WRI via NADMM without and with BM3D regularization. We perform FWI with the L-BFGS quasi-Newton method with line search to perform NADMM. 
We perform the inversion with noiseless data and set the maximum number of iteration to 70 as stopping criterion for IR-WRI in both of the cases (without and with BM3D). For FWI, the stopping criteria is set to the model error ($l_2$-norm of the difference between true and estimated model) achieved by IR-WRI for a fair comparison between the two waveform inversion methods.
Figure \ref{fig:Inclusions} shows the results obtained by FWI and IR-WRI for all four models without and with regularization.
It is clearly seen that for both methods regularization improved the results and successfully recovered different shapes of the anomalies, thanks to the adaptive nature of the BM3D. 
Although IR-WRI performed better than FWI, in this paper, we are not going to compare these methods because this is considered in \cite{Aghamiry_2019_IWR}, instead, we want to show how adaptive regularization can improve the results of these methods when it is applied using the proximal Newton algorithms. 

\begin{figure}[!htb]
\includegraphics[width=8cm,clip=true,trim=-0cm 0cm 0cm 0cm]{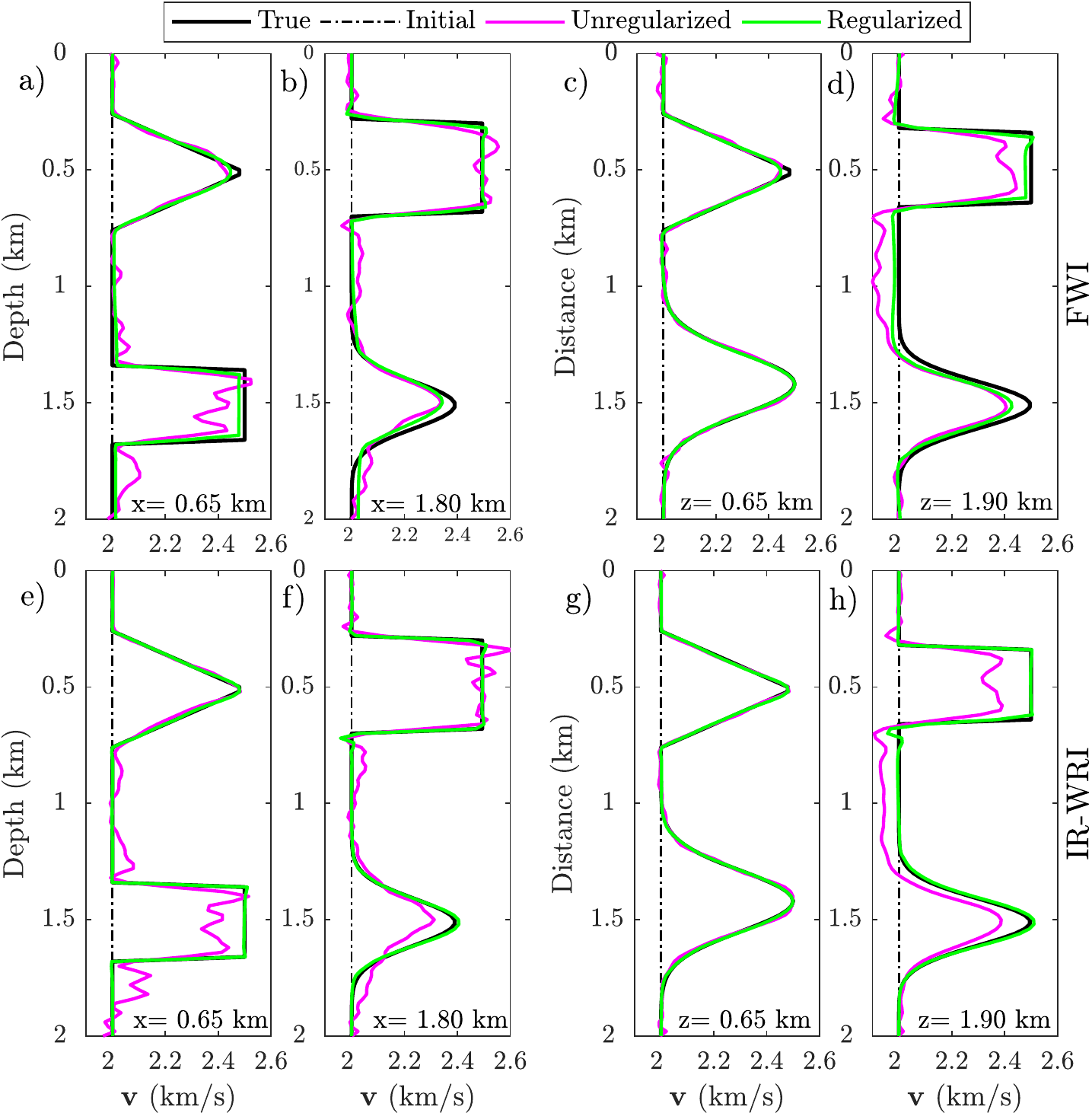}\\
\caption{Inclusion test. (a-b) Vertical logs at $x=0.65~km$ and $x=1.80~km$. (c-d) Horizontal logs at $z=0.65~km$ and $z=1.90~km$ for FWI results (true model is solid black, initial model is dashed black, the estimated model without regularization is blue and estimated model with regularization is red). (e-h) Same as (a-d) but for IR-WRI results.}
\label{fig:simple_log}
\end{figure}

We continue by using a model that includes all four inclusions (Fig. \ref{fig:simple}a). 
Figure \ref{fig:simple} shows the velocity models estimated by FWI (Figs. \ref{fig:simple}b-c) and IR-WRI (Figs. \ref{fig:simple}d-e) with and without regularization. 
A direct comparison between the true model, the initial model, and the final models without/with regularization along with two vertical logs at horizontal distances 0.65~km, 1.80~km, and two horizontal logs at vertical depths 0.65~km and 1.9~km are shown in Fig. \ref{fig:simple_log}.\\
\subsubsection{Robustness against noise}
\begin{figure}[!htb]
\centering
\includegraphics[width=0.3\textwidth]{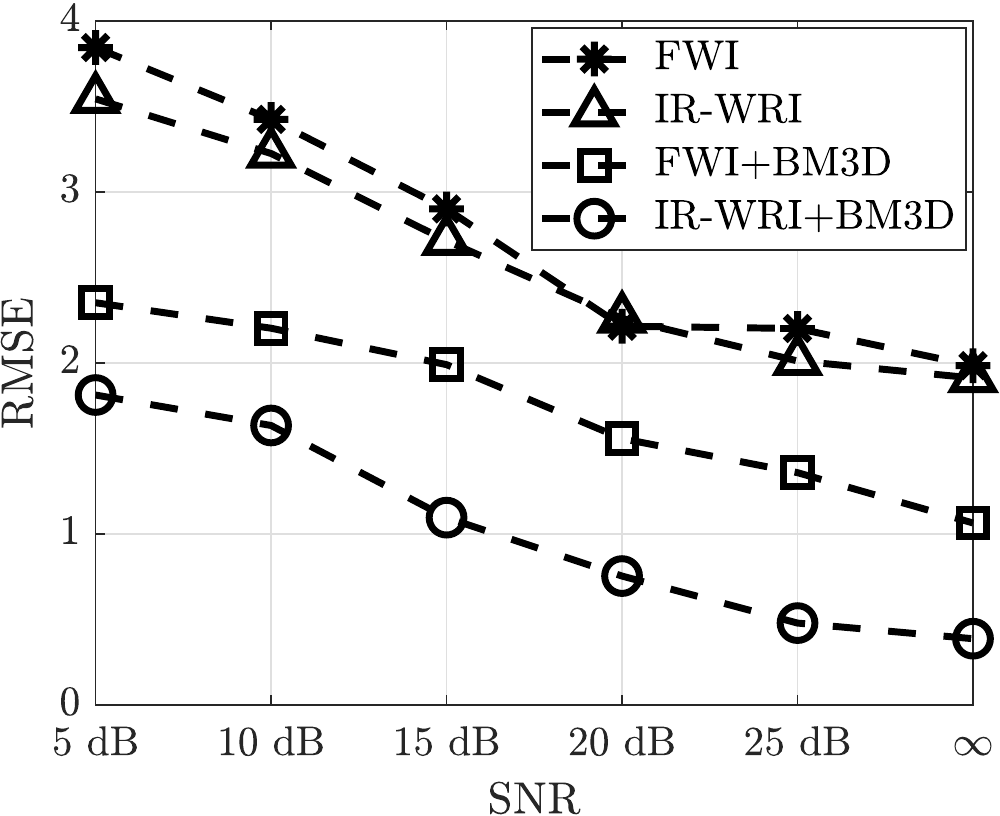}\\
\caption{RMSE for FWI and IR-WRI without and with BM3D using NADMM when data are contaminated with different level of noises.}
\label{fig:noise_sensitivity}
\end{figure}

\begin{figure} [!htb]
\includegraphics[width=0.49\textwidth]{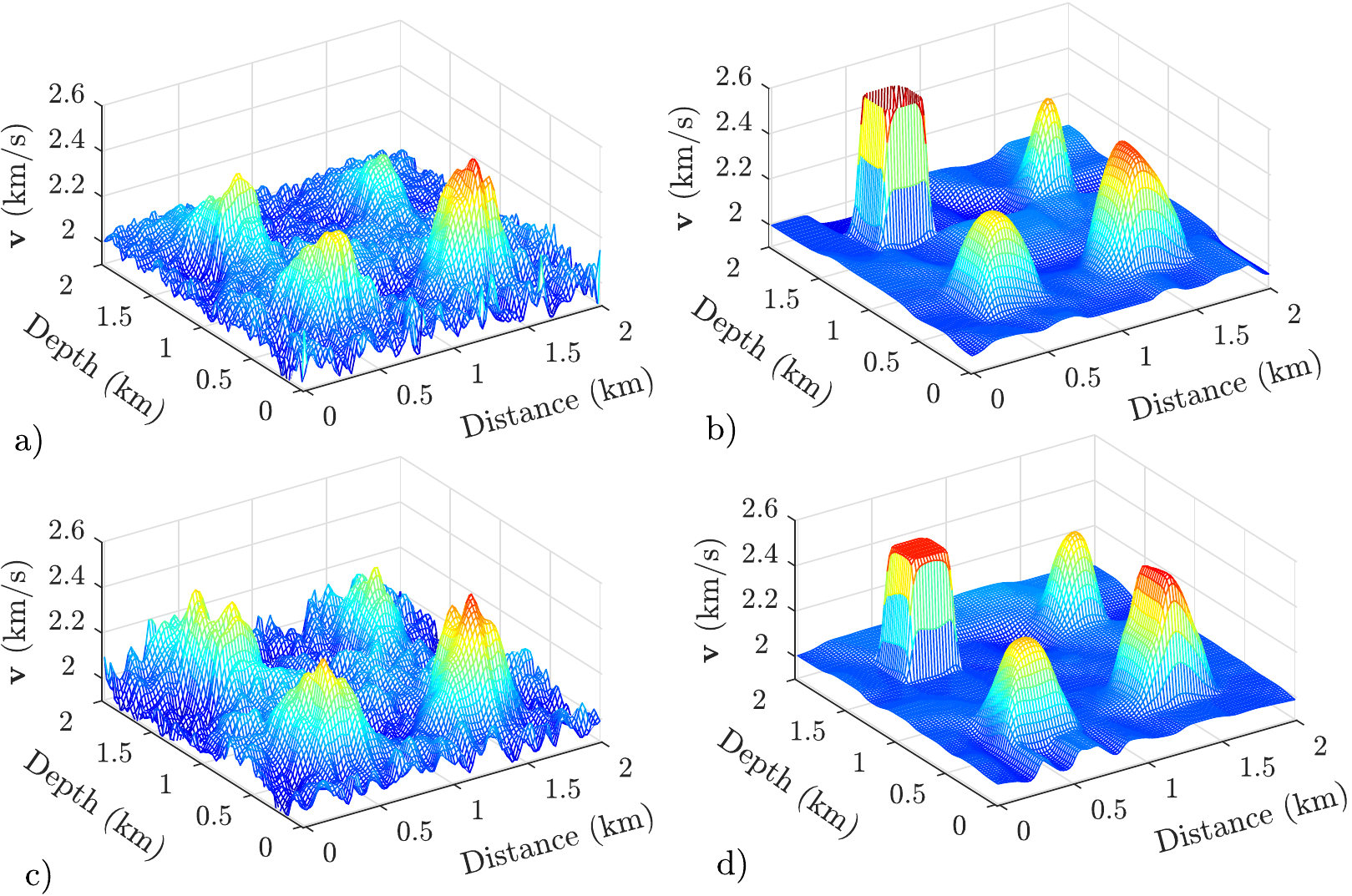}\\
\caption{The simple model test with SNR=$5$db. Estimated velocity model using (a) FWI without regularization, (b) FWI with regularization, (c) IR-WRI without regularization, (d) IR-WRI with regularization.}
\label{fig:simple_noisy}
\end{figure}

\begin{figure} [!htb]
\includegraphics[width=8cm,clip=true,trim=-0cm 0cm 0cm 0cm]{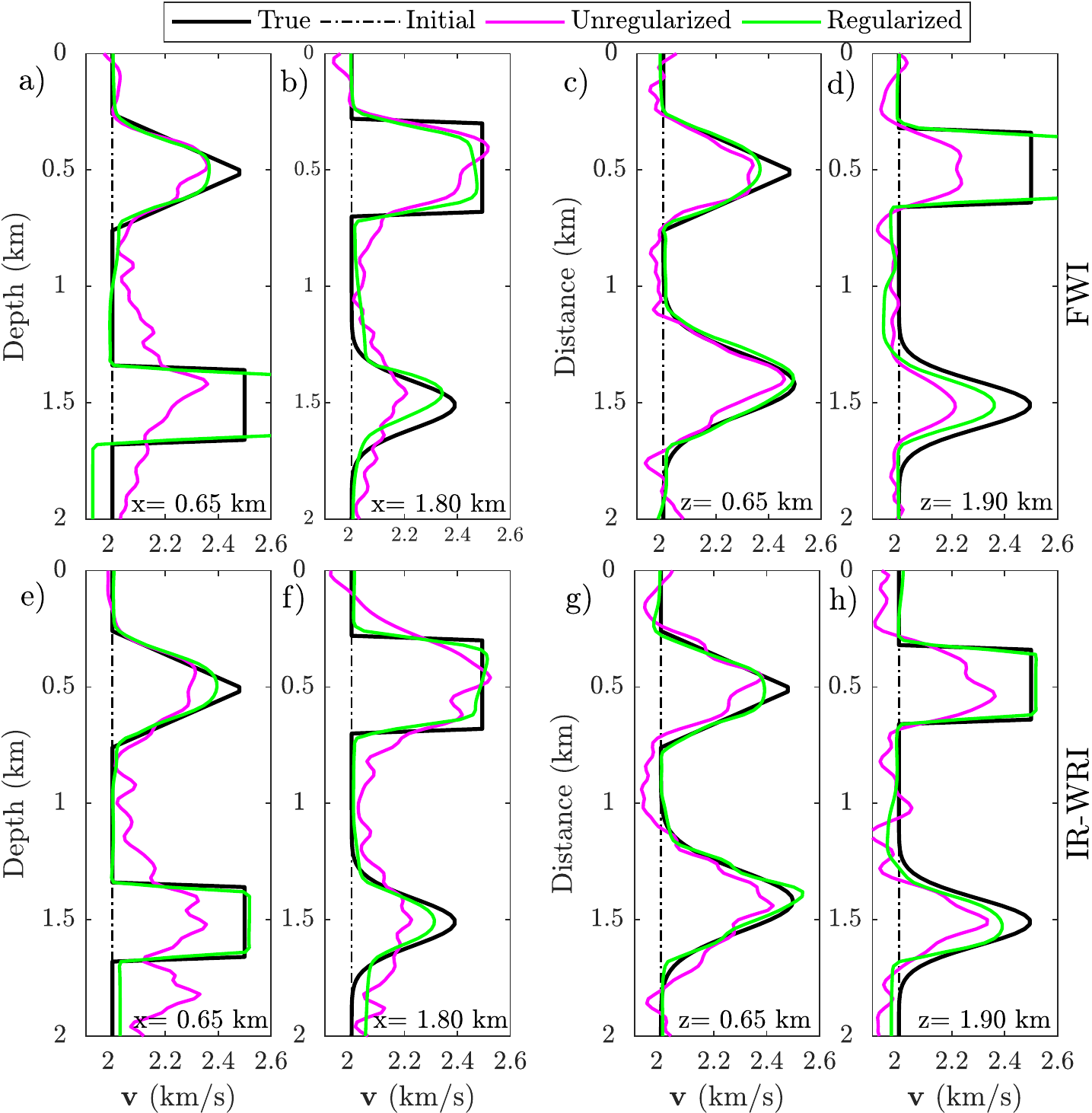}\\
\caption{Direct comparison of estimated velocity models with SNR=5db (Fig. \ref{fig:simple_noisy}). The configuration of this figure is the same as Fig. \ref{fig:simple_log}.}
\label{fig:simple_log_noisy}
\end{figure}

\begin{figure} [!htb]
\includegraphics[width=0.45\textwidth]{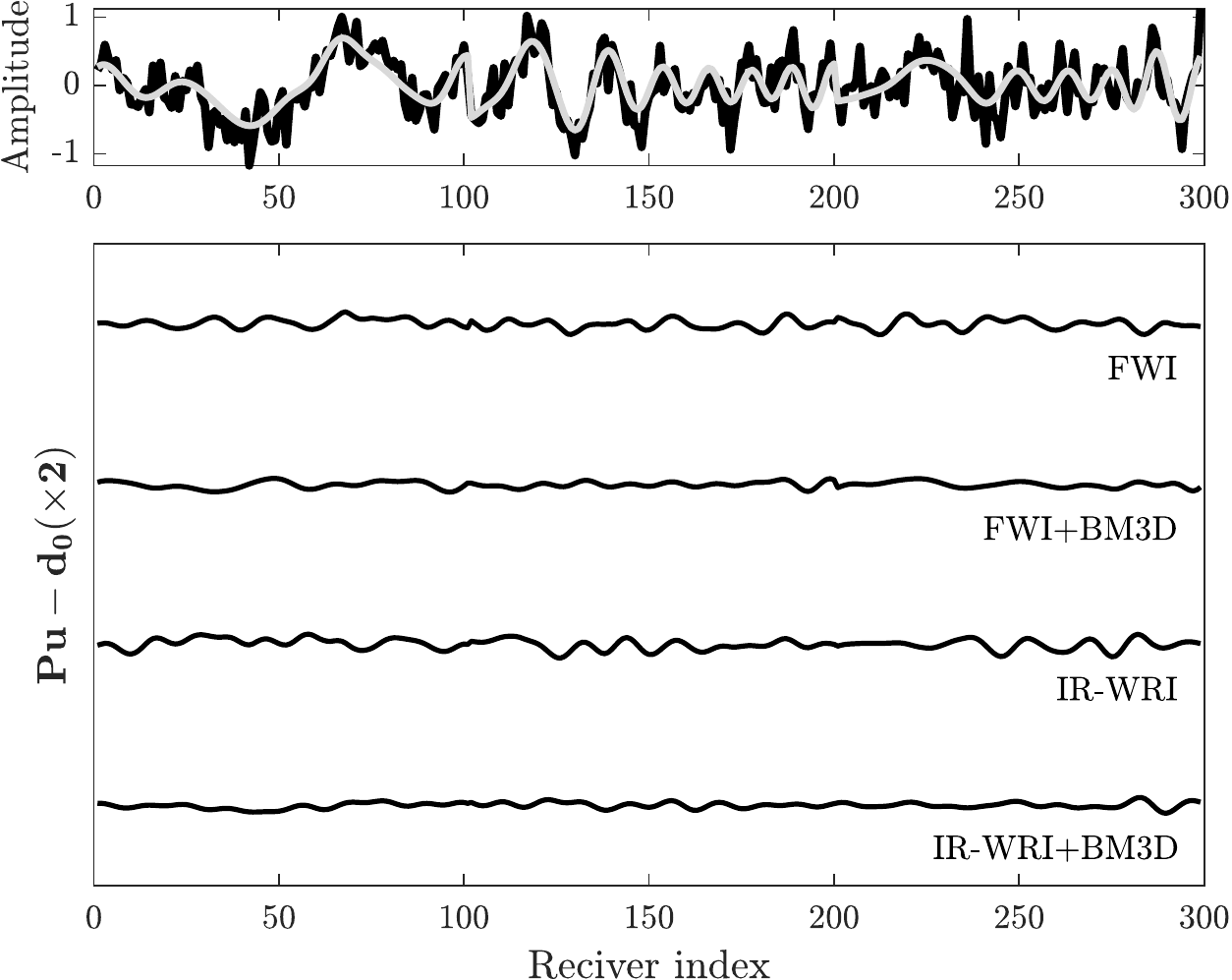}\\
\caption{Real part of the 10~Hz data. (a) The noisy data with SNR=5db are shown in black, while the noiseless data are shown in gray. (b) The difference between predicted data ($\bold{Pu}_k$) and noiseless data ($\bold{d}_0$) at the final iteration of Fig. \ref{fig:simple_noisy}. The residual curves are scaled by factor 2. }
\label{fig:data_fit}
\end{figure}

We continue by assessing the robustness of the proposed method against random noise. We apply FWI and IR-WRI without and with BM3D using NADMM when the data are contaminated with different level of random noises. The relative root mean square error (RMSE) curves versus signal to noise ration (SNR) is depicted in Fig. \ref{fig:noise_sensitivity}, where RMSE and SNR are defined as 
\begin{equation}
\text{RMSE}= 100\frac{\|\bold{m}-\bold{m}_*\|_2}{\|\bold{m}_*\|_2} ,
\end{equation}
in which $\bold{m}$ and $\bold{m}_*$ are the estimated and true models, respectively, and
\begin{equation}
\text{SNR}=20 \log \left( \frac{Signal~RMS~amplitude}{Noise~RMS~amplitude} \right).
\end{equation}
Fig. \ref{fig:noise_sensitivity} shows the average value (over 20 runs) for each SNR. Furthermore, we use $\|\bold{Pu}_k-\bold{d}\|_2 =1.01 \varepsilon$ as the stopping criterion of iteration, where $\varepsilon$ is the $\ell_2$ norm of the noise.
The velocity models estimated by FWI and IR-WRI without/with BM3D regularization for SNR=5db are shown in Fig. \ref{fig:simple_noisy}
with a direct comparison of the results in Fig. \ref{fig:simple_log_noisy}. 
In order to show how the data are fitted, the difference between the estimated and noiseless (10 Hz) data for different tests of Fig. \ref{fig:simple_noisy} are shown in Fig. \ref{fig:data_fit}.
\subsubsection{A comparison between NISTA and  NADMM}
\begin{figure} [!htb]
\includegraphics[width=8cm,clip=true,trim=-0cm 0cm 0cm 0cm]{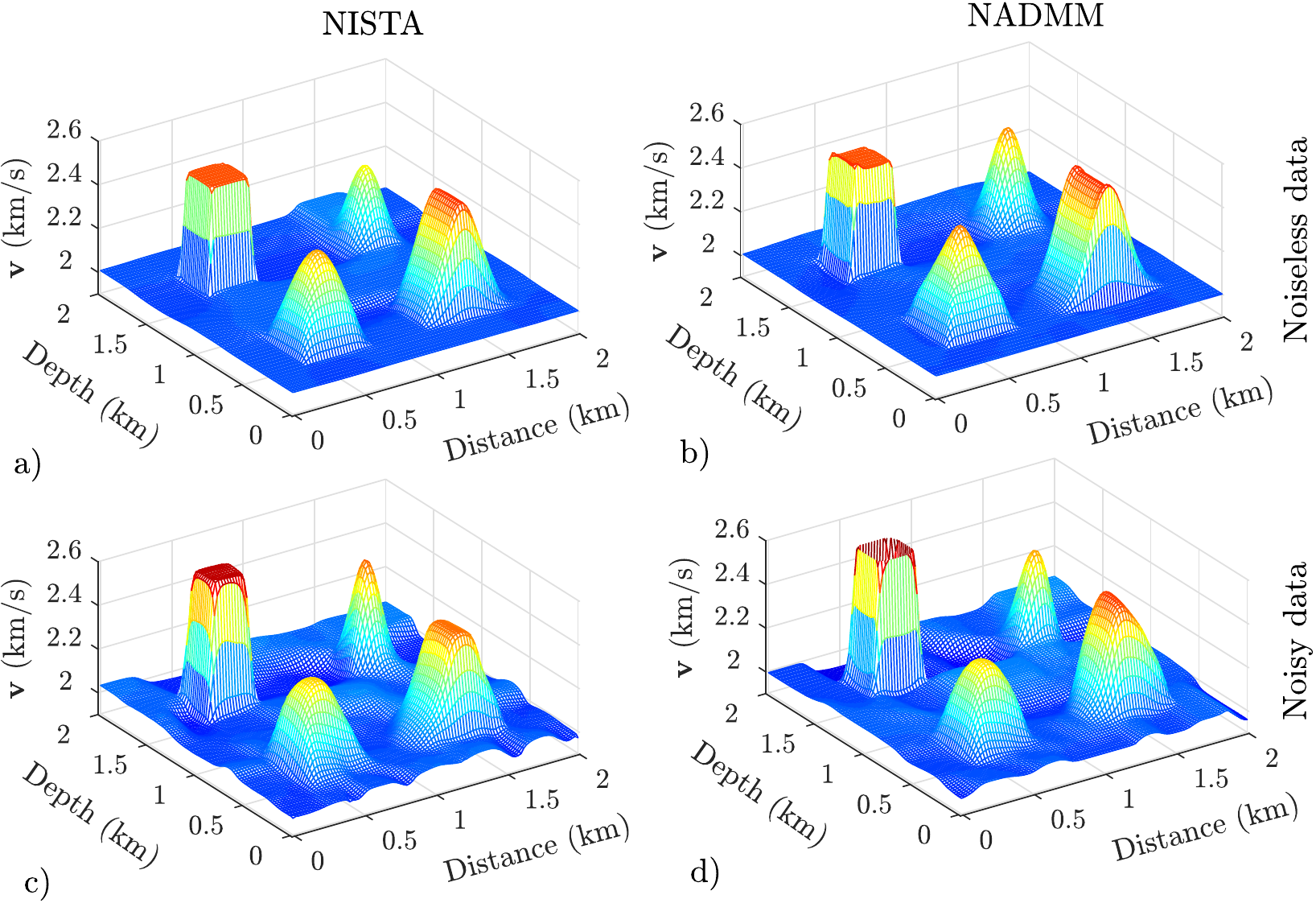}\\
\caption{BM3D regularized FWI on simple model using NISTA (a and c) and NADMM (b and d). (a-b) Noiseless data, (c-d) Noisy data with SNR=5db.}
\label{fig:simple_com_alg12}
\end{figure}

\begin{figure} [H]
\includegraphics[width=9cm,clip=true,trim=-0cm 0cm 0cm 0cm]{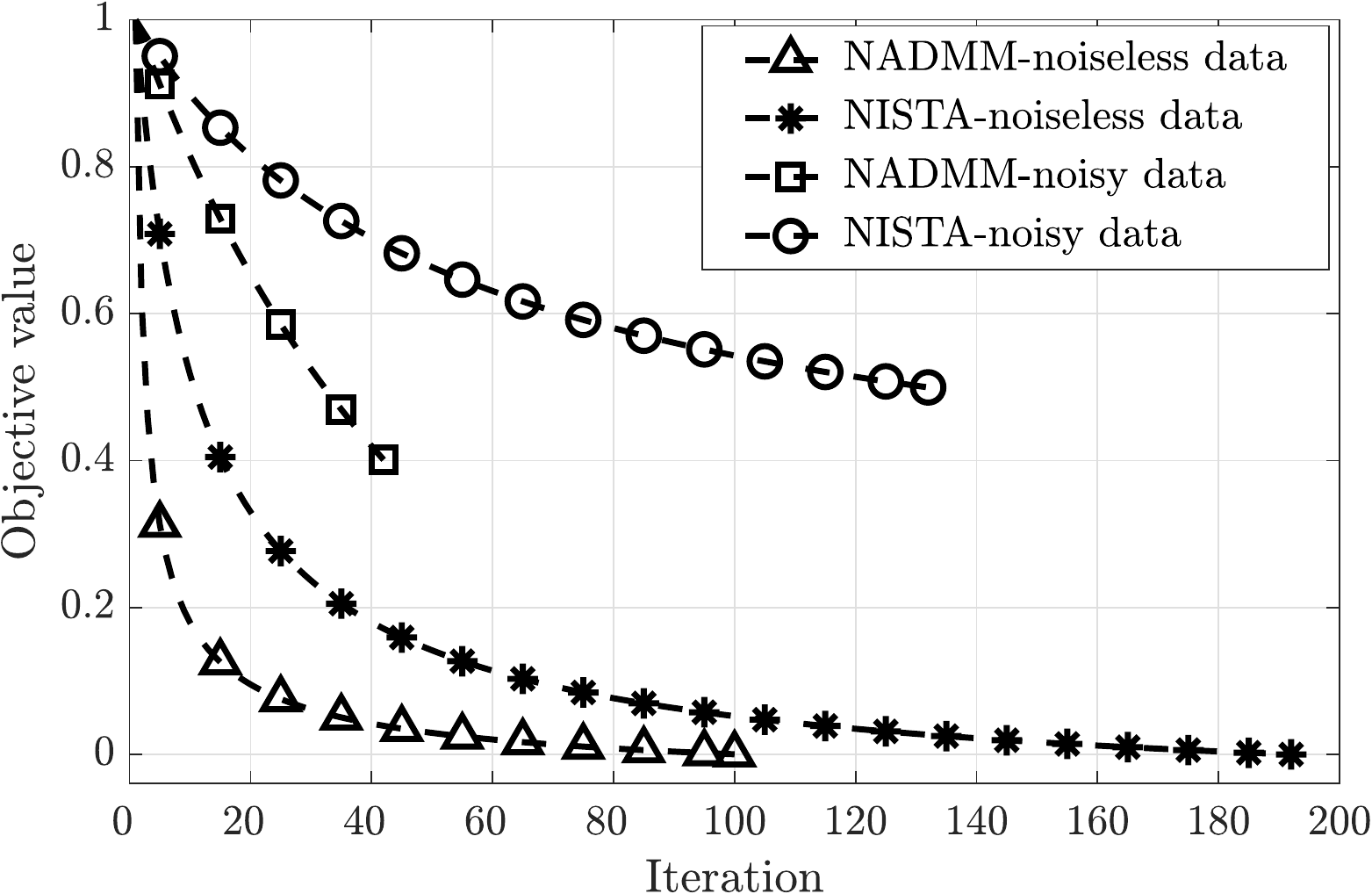}\\
\caption{Evaluation of the objective function of four different tests of Fig \ref{fig:simple_com_alg12}.}
\label{fig:simple_com_obj}
\end{figure}

Here we use the BM3D regularized FWI with noiseless and noisy data to compare NISTA and  NADMM. 
Fig. \ref{fig:simple_com_alg12} shows the estimated models obtained by both algorithms and Fig. \ref{fig:simple_com_obj} shows the corresponding convergence history (the objective function value) during the iterations.
Although we perform approximately 100 inner iterations of proximal gradient to estimate the search direction of NISTA, NADMM still performs better. Furthermore, since we implement both algorithms with L-BFGS, the results show that in practice NADMM should be preferred to NISTA.
%
%
%

%
%
%
%
%
\subsection{Performance on benchmark models}
\begin{figure*}[!htb] 
\includegraphics[width=1\textwidth]{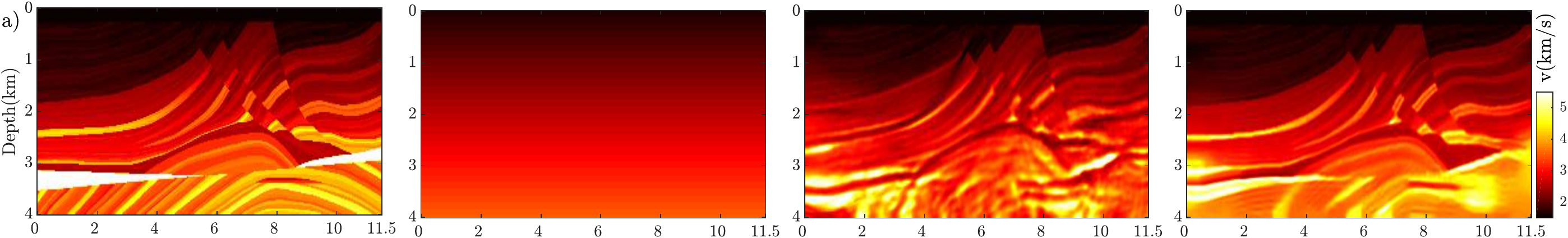}\\
\includegraphics[width=1\textwidth]{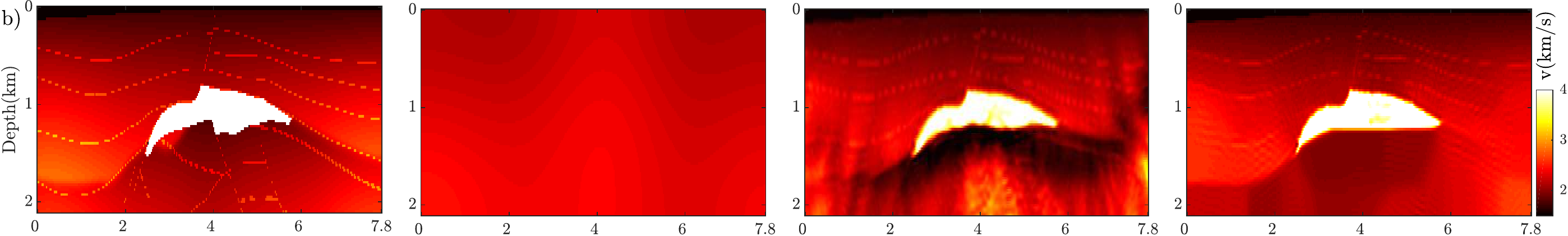}\\
\includegraphics[width=1\textwidth]{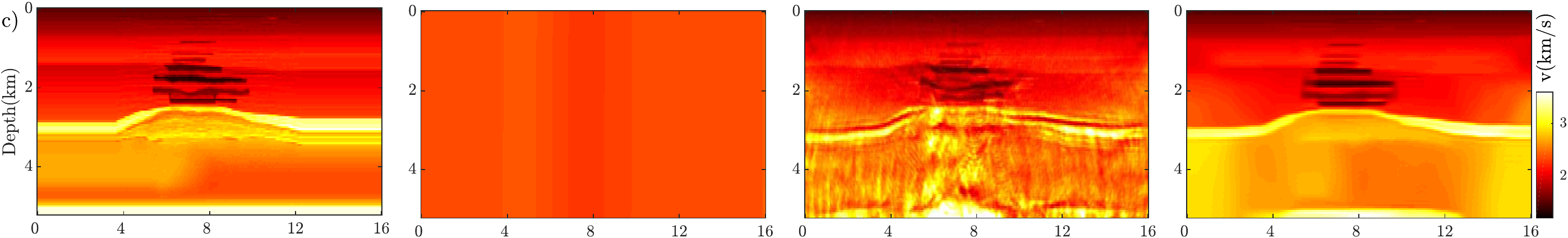}\\
\includegraphics[width=1\textwidth]{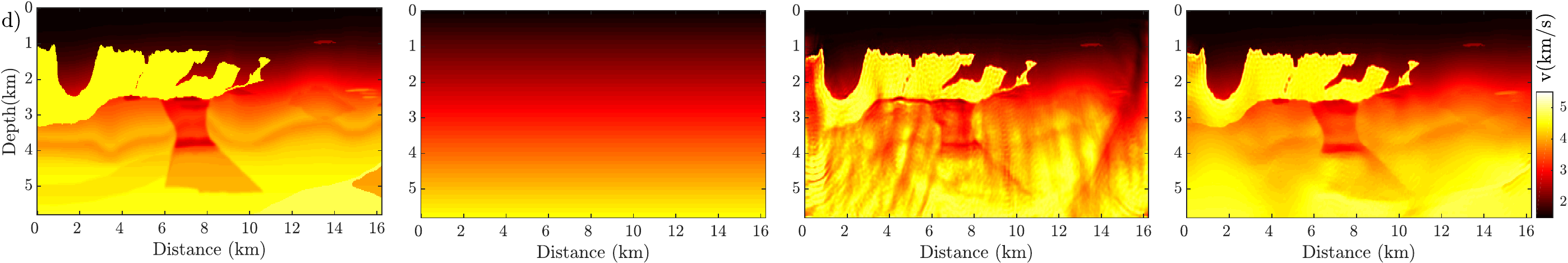}
\caption{IR-WRI without/with BM3D regularization on benchmark models. (a) Marmousi II, (b) SEG/EAGE salt, (c) Synthetic valhall and (d) 2004 BP salt models. The columns of this figure are as follows: True velocity model, initial velocity model, IR-WRI and BM3D regularized IR-WRI.}
\label{fig:bp_test}
\end{figure*}

\begin{figure*}[!htb] 
\includegraphics[width=0.48\textwidth]{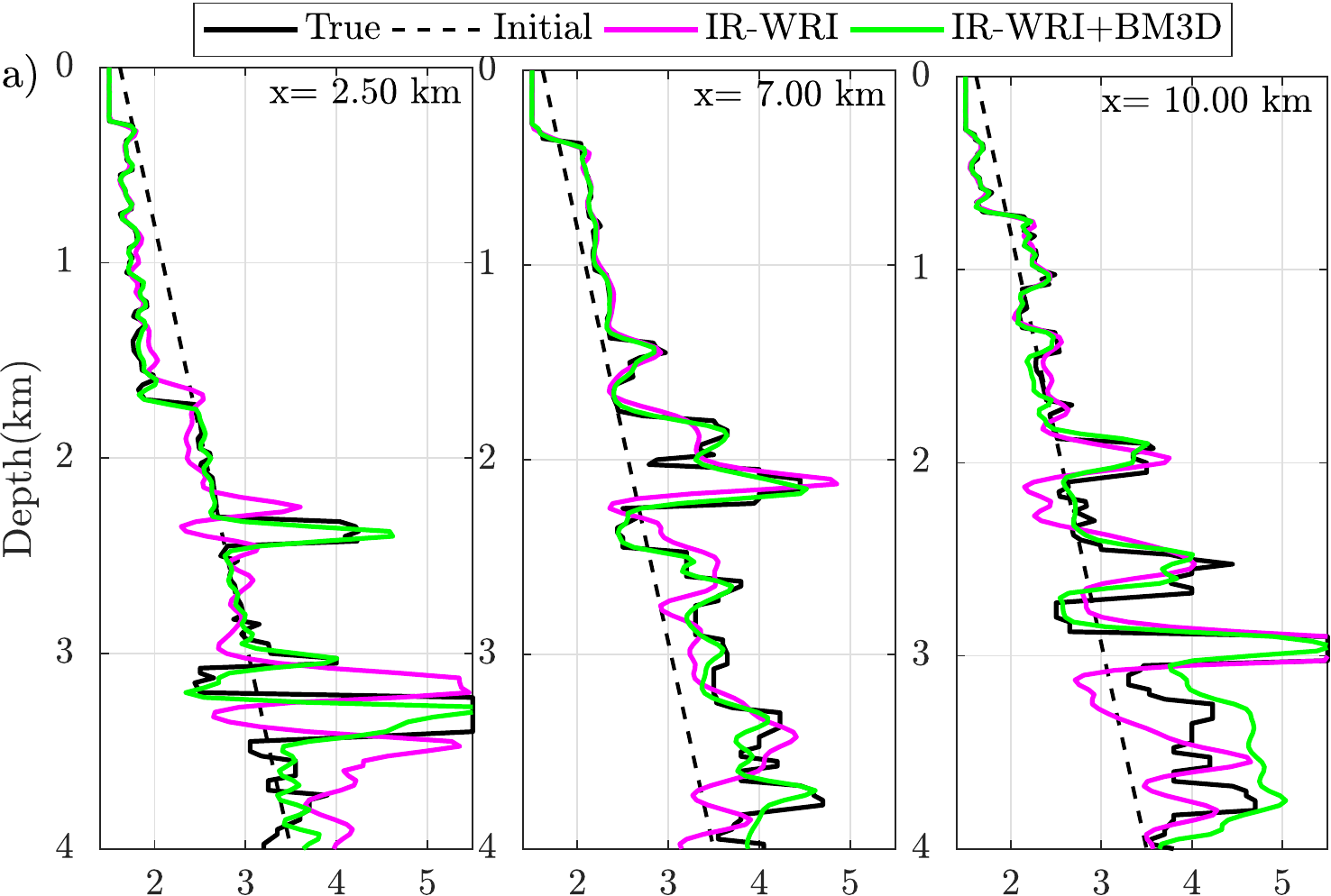}
\includegraphics[width=0.48\textwidth]{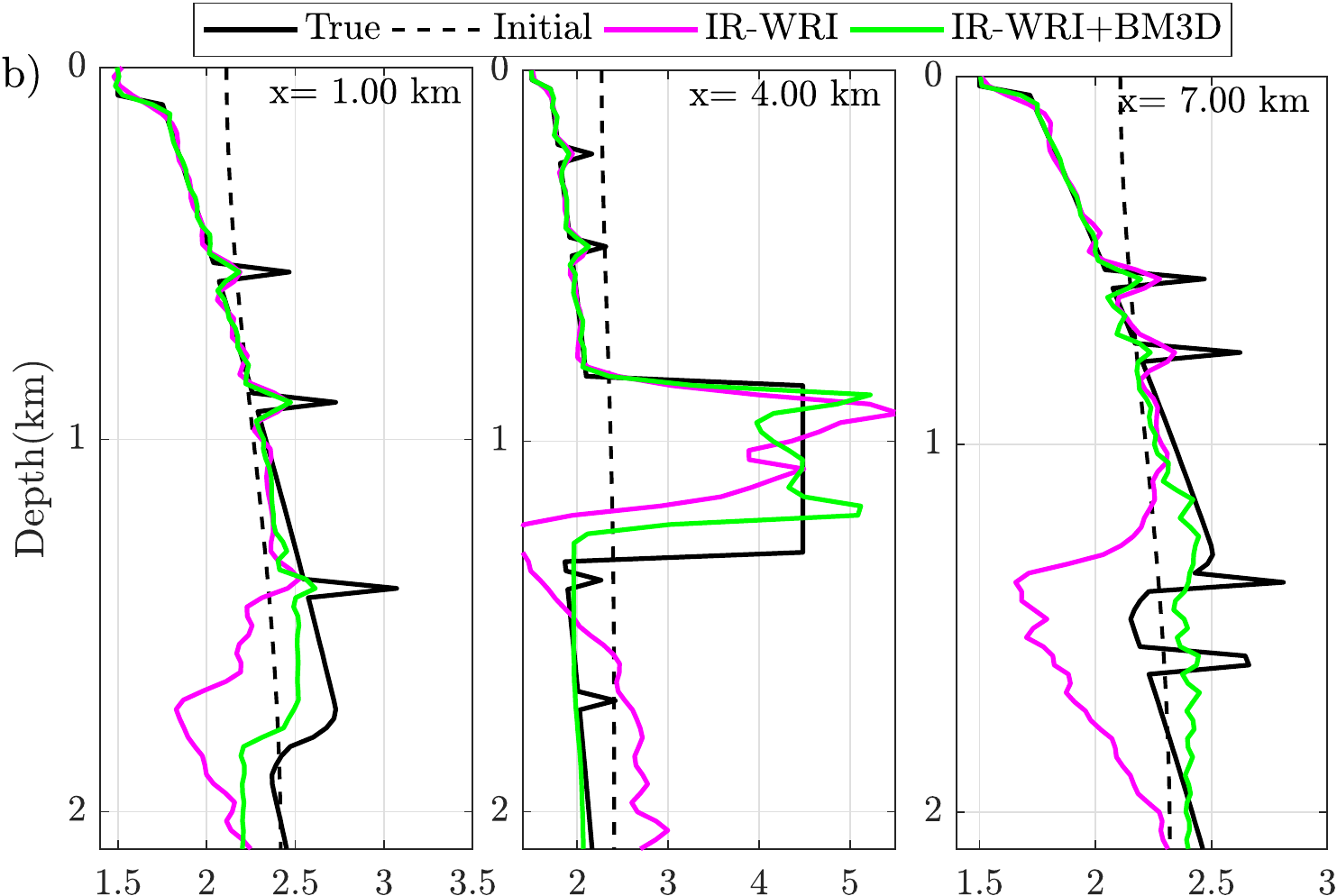}\\
\includegraphics[width=0.48\textwidth]{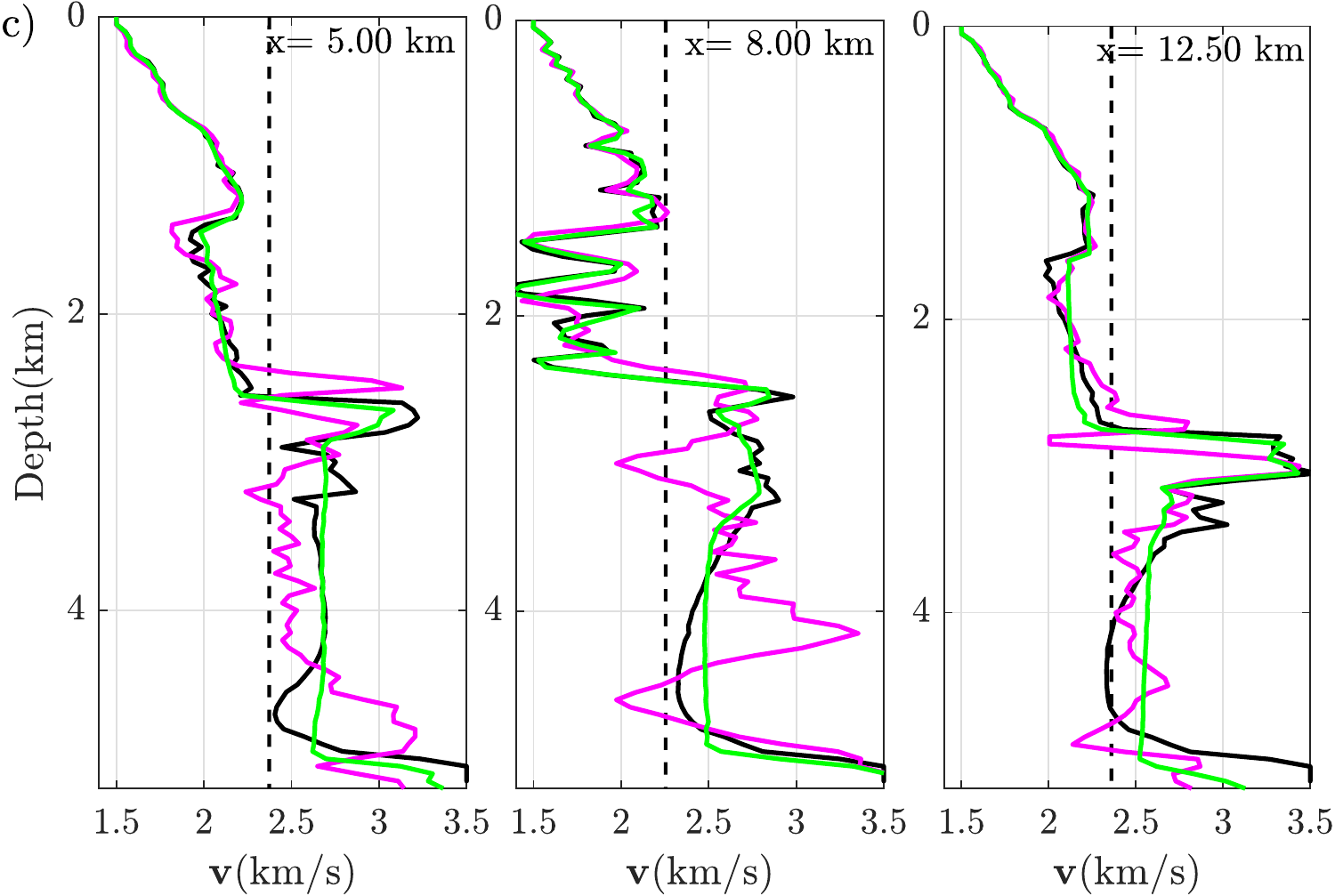}
\includegraphics[width=0.48\textwidth]{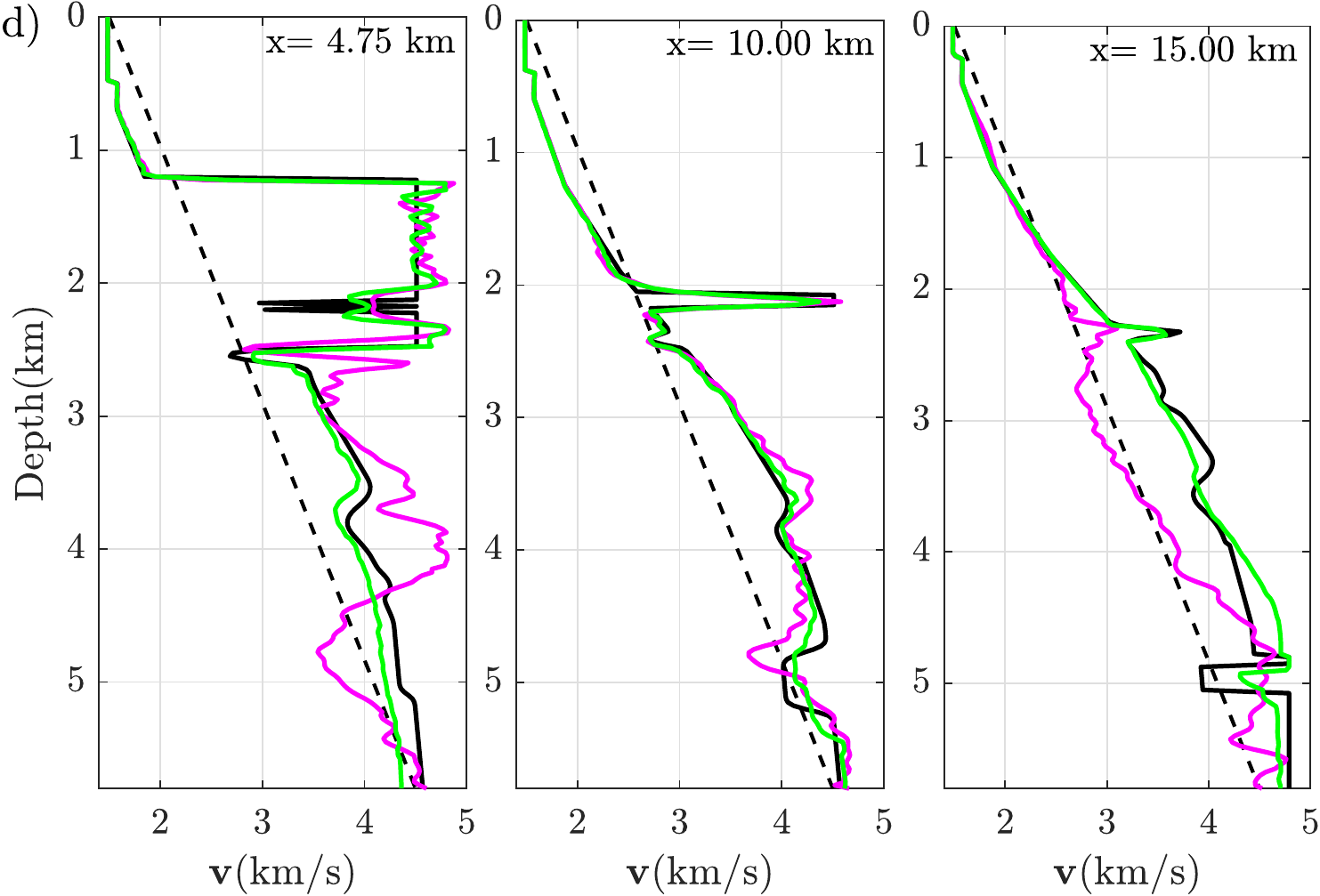}\\
\caption{Direct comparisons between velocity models estimated by IR-WRI without/with BM3D regularization . (a) Marmousi II, (b) SEG/EAGE salt, (c) Synthetic valhall and (d) 2004 BP salt models. In each panel, the true model is solid black, initial model is dashed black, the estimated model without regularization is pink and estimated model with regularization is green. The horizontal location of each log is written in each panel.}
\label{fig:bp_test_log}
\end{figure*}

We continue by assessing the performance of the proposed adaptive IR-WRI algorithm using more complicated models when the models are selected from well-documented 2D benchmark subsurface velocity models in exploration seismic, e.g. the Marmousi II \cite{Martin_2006_M2E}, SEG/EAGE overthrust \cite{Aminzadeh_1997_DSO}, SEG/EAGE salt \cite{Aminzadeh_1997_DSO}, synthetic Valhall \cite{Prieux_2011_FAI} and 2004 BP salt \cite{Billette_2004_BPB} benchmark velocity models. The selected target from these benchmark models are shown in the first column of Fig. \ref{fig:bp_test}, respectively.
The fixed-spread acquisition with a few equally spaced sources at the sea bottom and a line of equally spaced receivers at the depth 25~m is used for all of the tests. Also, the models are discretized with 25 m spacing in horizontal and vertical directions (see Table \ref{tab:info} for more technical details). We compute the wavefields using Perfectly-Matched Layer (PML) absorbing boundary conditions along the bottom, right, and left sides of the model using 10 grid points in the PMLs and a free-surface boundary condition at the surface when a 10~Hz Ricker wavelet is used as the source signature. We design a multiscale inversion with a classical continuation frequency strategy in the selected frequency band by proceeding over small batches of two frequencies with a frequency interval of 0.5 Hz. We also perform three paths through the batches, where the starting and finishing frequencies of the paths and other technical details about the modes are reported in Table \ref{tab:info}. The initial velocity models are crude models, as shown in the second column of Fig. \ref{fig:bp_test}. Accordingly, we tackle these benchmarks with IR-WRI only since FWI would remain stuck in a local minimum due to cycle skipping. We set the number of IR-WRI iterations per frequency batch equal to 10 or $\ell_2$-norm of source residuals equal to 1e-3 as the stopping criteria. %
The estimated models without and with BM3D regularization are shown in the third and fourth columns of Fig. \ref{fig:bp_test}, respectively.
A direct comparison between the true velocity, the initial and the final velocity models without/with BM3D regularization are shown in Fig \ref{fig:bp_test_log}a-d, for Marmousi II, SEG/EAGE salt, Synthetic valhall and 2004 BP salt models, respectively. The results show that, although the different benchmark models are characterized by different kinds of structures, adaptive regularization combined with IR-WRI manages to reconstruct accurately each of them with a significant jump of quality compared to the case where BM3D is not used.\\ 

\begin{table*}
\caption{Technical details of the benchmark models}
\label{tab:info}
\begin{tabular}{c|c|c|c|c|c|}
\cline{2-6}
                                          & \begin{tabular}[c]{@{}c@{}}Size \\ (km $\times$ km)\end{tabular} & \begin{tabular}[c]{@{}c@{}}Source \\ interval(m)\end{tabular} & \begin{tabular}[c]{@{}c@{}}Receiver\\  interval(m)\end{tabular} & \begin{tabular}[c]{@{}c@{}}Inverted frequency band\\ {[}starting-finishing{]} frequencies(Hz)\end{tabular} & \begin{tabular}[c]{@{}c@{}}Outer iterations\\ {[}starting-finishing{]} frequencies(Hz)\end{tabular} \\ \hline
\multicolumn{1}{|c|}{Marmousi II}         & 4.25 $\times$ 11.5                                               & 250                                                           & 50                                                             & {[}3-10{]}                                                                                             & {[}3-8{]}, {[}4-9{]}, {[}5-10{]}                                                                 \\ \hline
\multicolumn{1}{|c|}{SEG/EAGE salt model} & 2.1 $\times$ 7.8                                                 & 100                                                           & 25                                                             & {[}3-7{]}                                                                                              & {[}3-6{]}, {[}3.5-7{]}, {[}4-7{]}                                                                \\ \hline
\multicolumn{1}{|c|}{Synthetic valhall}   & 5.25 $\times$ 16                                                 & 500                                                           & 100                                                            & {[}3-13{]}                                                                                             & {[}3-9{]}, {[}4-11.5{]}, {[}5-13{]}                                                              \\ \hline
\multicolumn{1}{|c|}{2004 BP salt}        & 5.8 $\times$ 16.25                                               & 250                                                           & 50                                                             & {[}3-13{]}                                                                                             & {[}3-9.5{]}, {[}3.5-11.5{]}, {[}5-13{]}                                                          \\ \hline
\end{tabular}
\end{table*}
\section{Conclusions} 
In this paper, we proposed a flexible framework to apply state-of-the-art regularizations embedded in denoising algorithms to nonlinear inverse problems, in particular full-waveform inversion and its variant by wavefield reconstruction. In these proximal Newton-type algorithms, the the search direction implicitly involves the gradients/subgradients of the possibly non-differentiable regularization function.
The regularization is treated as black box, allowing sophisticated regularizers to be employed via denoising engines.  
Two proximal Newton algorithm have been proposed. The first relies on FISTA, while the second relies on ADMM. Numerical tests with the Rosenbrock function and a toy example with multiple inclusions of different shape suggest that ADMM provides the fastest convergence.
Several numerical examples using wavefield inversion with L-BFGS and BM3D denoiser were tested confirming that the proposed proximal Newton-type algorithms successfully recovers complicated velocity models without needing prior information about the targeted structure. 

\section{Acknowledgments}  
This study was partially funded by the WIND consortium (\textit{https://www.geoazur.fr/WIND}), sponsored by Chevron, Shell, and Total, as well as by the IDEX UCA JEDI with the WIMAG project. This study was granted access to the HPC resources of SIGAMM infrastructure (http://crimson.oca.eu), hosted by Observatoire de la C\^ote d'Azur and which is supported by the Provence-Alpes C\^ote d'Azur region, and the HPC resources of CINES/IDRIS/TGCC under the allocation A0050410596 made by GENCI."
\IEEEpeerreviewmaketitle

\bibliographystyle{IEEEtran}
\bibliography{IEEEabrv}
\newcommand{\SortNoop}[1]{}

\end{document}